\documentclass[a4paper, 12pt, oneside, article]{memoir}
\usepackage{palgrave}
\bibliography{Set.bib}
\title{The iterative notion of function and the iterative notion of set}
\author{Tim Button}
\emailaddress{tim.button@ucl.ac.uk}
\begin{document}\midsloppy
\pagestyle{nicelypage}
\maketitletop\selectbodyfont\noindent
\textcolor{blue}{This is a preprint and may be subject to minor alterations. Forthcoming in  C.\ Antos, N.\ Barton, and G.\ Venturi (eds.), \emph{Palgrave Companion to the Philosophy of Set Theory}.}

\
\\
Hilary Putnam once suggested that ``the actual existence of sets as `intangible objects' suffers\ldots from a generalization of a problem first pointed out by Paul Benacerraf\ldots are sets a kind of function or are functions a sort of set?''\footnote{\textcite[13.Dec.2014]{Putnam:blog}.} Sadly, he did not elaborate; my aim here is to do so on his behalf. 

There are well-known methods for treating sets as functions and functions as sets. But these do not raise any obvious philosophical or foundational puzzles (see \S\ref{s:vonNeumann}). To raise such puzzles, we first need to provide a full-fledged function theory. I will supply one, which axiomatizes the iterative notion of function in exactly the same sense that \ZF axiomatizes the iterative notion of set. Indeed, this theory is synonymous with \ZF (see \S\S\ref{s:story}--\ref{s:flt}).

It might seem that set theory and function theory present us with rival foundations for mathematics, since they postulate different ontologies. But I will argue that appearances are deceptive. Set theory and function theory provide the very same judicial foundation for mathematics (see \S\ref{s:judicial}). They do not supply rival metaphysical foundations (see \S\ref{s:metaphysics}); indeed, if they supply metaphysical foundations at all, then they supply the very same metaphysical foundations (see \S\ref{s:metaphysicsnotdead}).

\section{Sets as functions / functions as sets}\label{s:vonNeumann}
There is a standard way to regard functions as sets. To avoid ambiguity, say that a set $f$ is a \emph{\setfunction} iff both (i) every member of $f$ is a Kuratowski ordered-pair, and (ii) if $\tuple{x,y} \in f$ and $\tuple{x,z} \in f$ then $y=z$. We then treat ``$f(x) = y$'' as an abbreviation of ``$\tuple{x,y} \in f$''. Conversely, sets simply can be regarded as \emph{characteristic} functions, that is, as functions whose values are always either $1$ or $0$. We then treat ``$x \in y$'' as an abbreviation for ``$y(x) = 1$''. 

This simple observation suggests a formal equivalence. Consider the following re-axiomatization of \ZF. Replace each atomic formula ``$x \in y$'' with ``$\appof{y}{x} = 1$'', where ``$\appfunction$'' is a new two-place function symbol. We can read ``$\appof{y}{x}$'' as ``$y$ applied to $x$''. Next, add an axiom:
	$$\forall y \forall x(\appof{y}{x} = 0 \lor \appof{y}{x} = 1)$$
Finally, treat $0$ and $1$ as objects such that: $\appof{0}{x} = 0$ for all $x$; and $\appof{1}{x} = 1$ iff $x = 0$. Call the result \ZFchi. Fairly trivially, \ZFchi is \emph{synonymous} with \ZF. 

I should unpack this notion of synonymy. Two theories, \Thvar{S} and \Thvar{T}, are said to be \emph{synonymous} when each interprets the other, and running the two interpretations consecutively gets you back where you began. More precisely, it is to say that there are translations, $\intTtoS$ and $\intStoT$, which are \emph{interpretations}, i.e.\ theorems are translated to theorems, so that:
\begin{listn-0}
	\item\label{dag:int} if $\Thvar{T} \proves \phi$ then $\Thvar{S} \proves \phi^\intTtoS$, for any 
	sentence $\phi$ in \Thvar{T}'s language; and
	\item\label{ddag:int} if $\Thvar{S} \proves \phi$ then $\Thvar{T} \proves \phi^\intStoT$, for any sentence $\phi$ in \Thvar{S}'s language
\end{listn-0}
Moreover, composing these interpretations gets you back where you started, in the sense that:\footnote{I write $\phi^{IJ}$ for $(\phi^I)^J$. Note that $\phi$ and $\phi^{\intStoT\intTtoS}$ may be syntactically distinct, but they are provably equivalent. For more, see e.g.\ \textcites[\S\S1--2]{FriedmanVisser:WBIS}[ch.5]{ButtonWalsh:PMT}.}
\begin{listn}
	\item\label{ddag:dag} $\Thvar{S} \proves \phi \liff \phi^{\intStoT\intTtoS}$, for any formula $\phi$ in \Thvar{S}'s language; and
	\item\label{dag:ddag} $\Thvar{T} \proves \phi \liff \phi^{\intTtoS\intStoT}$, for any formula $\phi$ in \Thvar{T}'s language
\end{listn}

Now, \ZF is a theory of sets, and \ZFchi is a theory of (characteristic) functions; but \ZF and \ZFchi are \emph{synonymous}, in the sense just defined. And something like this observation underwrote John von Neumann's extremely suggestive claim that the concepts \emph{set} and \emph{function} ``are completely equivalent, since a function can be regarded as a set of pairs, and a set as a function that can take two values.''\footnote{\label{fn:vNimplementation}\textcite[221--2/396]{vonNeumann:EAM}; see also \parencites*[230--1/411--12]{vonNeumann:EAM}[676]{vonNeumann:DAM}[227, 231]{vonNeumann:WAM}.}

That conclusion is too quick. For one thing: there is nothing especially \emph{set-theoretic} about that equivalence. Not just membership, but any relation, can be given this treatment: just read ``$R(x_1, \ldots, x_n, y)$'' as an abbreviation of ``$@_R(y, x_1, \ldots, x_n) = 1$'', and add an axiom:
	$$\forall y \forall x_1 \ldots \forall x_n(@_R(y, x_1, \ldots, x_n) = 0 \lor @_R(y, x_1, \ldots, x_n) = 1)$$
This provides a perfectly general method for rewriting a theory in functional terms. But it is needlessly baroque: it just amounts to incorporating surrogates for truth-values, $0$ and $1$, into the theory's object language.

For another thing, the synonymy of \ZF and \ZFchi would only suggest that the concept \emph{set} is intimately connected with the concept \emph{characteristic function}, not with \emph{function} simpliciter. Now, to be fair to von Neumann, he was not considering the theory \ZFchi, but a much richer theory which encompassed functions which are \emph{not} characteristic functions.\footnote{So, rather than having an \emph{axiom} like $\forall y \forall x(\appof{y}{x} = 0 \lor \appof{y}{x} = 1)$, he \parencites*[229/403]{vonNeumann:EAM}[678, 684]{vonNeumann:DAM} stipulatively defines that the \emph{sets} will be those (I-objects) $y$ such that $\forall x(\appof{y}{x} = 0 \lor \appof{y}{x} = 1)$. In \S2.2.1 of chapter 6, Toby Meadows presents a theory, vNC, which is similar to von Neumann's. (I first encountered Meadows' work in March 2023, after writing this paper.)} 
However, he also did not advance a synonymy result (the notion of \emph{synonymy} had not yet been invented).\footnote{That first occurs with \textcites[\S6]{Montague:PhD}{Bouvere:ST}.} So I suspect that von Neumann's thought was roughly that the concepts \emph{set} and \emph{function} are equivalent, because his rich function theory can interpret (an unspecified)\footnote{The theory \Thnamed{NBG} would do.} set theory via characteristic functions, and (that same) set theory can interpret his function theory via functions$_\in$. (So we would have conditions \eqref{dag:int} and \eqref{ddag:int} from above, but we would lack conditions \eqref{ddag:dag} and \eqref{dag:ddag}.) But if that was indeed von Neumann's thought, it falls short of what he needed. Intuitively, (mere) mutual interpretability is insufficient to establish the equivalence of the relevant concepts; after all, \PA and $\PA + \lnot\textspaced{Con}(\PA)$ are mutually interpretable, but they surely articulate very different conceptions of \emph{number}. 

All told, we do not \emph{yet} have a good reason to agree with von Neumann's conclusion that the concepts \emph{set} and \emph{function} are ``completely equivalent''. Nonetheless, von Neumann was onto something; there \emph{is} an important equivalence between these concepts. My aim in this paper is to unpack this equivalence. 

\section{The iterative notion of \emph{function}}\label{s:story}
I must start by making my target concept of \emph{function} more precise. In particular, and analogously with the iterative notion of \emph{set}, I will consider the iterative notion of \emph{function}.

Recall that we can introduce the iterative notion of \emph{set} with the following story:\footnote{This formulation is from \textcite{Button:LT1}; but such stories are hardly novel to me! In that paper, I explain connections to e.g.\ \textcites{Boolos:ICS}{Scott:AST}[323]{Shoenfield:AST}. Note that we need no more (nor less) than this Story to introduce the (bare) idea of the iterative notion of \emph{set}.}
\begin{storytime}\textbf{The Set Story.}	
	Sets are arranged in stages. Every set is found at some stage. At any stage $\stage{s}$: for any sets found before $\stage{s}$, we find a set whose members are exactly those sets. We find nothing else at $\stage{s}$.
\end{storytime}\noindent
This Story articulates the bare-bones iterative notion of \emph{set}. Adopting that conception of \emph{set} is, of course, a choice.\footnote{We might instead have pursued Quine's \NF.} But it is a common choice, and it has much to recommend it. I will say a bit more about its commendable qualities in \S\ref{s:judicialfoundations}; for now, I note that the iterative conception immediately blocks the threat of Russell's paradox.

Raphael Robinson noted that we can introduce an iterative notion of \emph{function} with a similar story:\footnote{``Suppose that we have the idea `function,' but nothing to use as arguments and values. We can start with the function $0$ which is not defined for any arguments and hence does not need any values\ldots. Then we construct a function $1$ such that $1(0) = 0$, but $1(x)$ is undefined for $x\neq 0$\ldots We proceed, always using as arguments and values functions already constructed.'' \parencite[29--30, notation slightly changed]{Robinson:TCMVNS}. Robinson then offered a deliberately simplified version of von Neumann's theory (which is quite different from my own \FLTzf). I wrote this paper in embarrassing ignorance of Robinson's paper; many thanks to Chris Menzel for bringing it to my attention.} 
\begin{storytime}\textbf{The Function Story.}
	Functions are arranged in {stages}. Every function is found at some stage. At any stage $\stage{s}$: we find all functions whose domains and ranges are exhausted by the functions we found earlier. We find nothing else at stage $\stage{s}$.
\end{storytime}\noindent 
Again, this Story lays down nothing more than the bare idea of an iterative notion of \emph{function}. Furthermore, embracing the Story, and thereby adopting the iterative notion of \emph{function}, is a choice.\footnote{We could, instead, have developed a notion of \emph{function} along the lines of Church's untyped $\lambda$-calculus; or along the lines of \citeauthor{Santanna:Flow}'s 
\parencite*{Santanna:Flow} theory Flow.} But it is a good choice. Apart from anything else, it immediately blocks the threat posed by the paradoxes of na\"ive function theory. To illustrate, let $d$ be this na\"ive diagonal function:\footnote{cf.\ \textcites[cf.][676--7 fn5]{vonNeumann:DAM}[31]{Robinson:TCMVNS}.}
\begin{align*}
	\text{if }x(x) \neq x&\text{, then }d(x) = x\\
	\text{if }x(x) = x&\text{, then }d(x) \text{ is undefined}
\end{align*}
Contradiction follows swiftly: if $d(d) = d$ then $d(d)$ is undefined, which is absurd; so $d(d) \neq d$, and hence $d(d)=d$, a contradiction. But, given the Function Story, $d$ obviously does not exist; its mapping behaviour requires it to be found after every function, which is absurd.

Two clarificatory comments might be helpful. The functions introduced in the Function Story are always \emph{partial} functions: for example, $f(f)$ is undefined for every $f$. Moreover, the functions are \emph{pure}: they only map from functions to functions. We could modify the Story to allow ``urelements''---and ultimately we probably should---but for simplicity, and for parity of presentation with the case of pure sets, I will ignore them in this paper.\footnote{\textcite[\S\S{}A--B]{Button:LT1} tackles the Set Story with urelements; the same techniques can be carried over to the Function Story.}

\section{Functional Level Theory}\label{s:flt}
Of course, working within \ZF, we can recursively define a system of (hereditary-)functions$_\in$ which would satisfy the Function Story. But I cannot rest content with what we can do within \ZF. After all, I want to show that set-theoretic and function-theoretic approaches are on a par when it comes to mathematical foundations. To do this, I need to axiomatize the Function Story without relying upon some ``prior'' theory of sets (like \ZF); I need to give you an axiomatic theory of functions which stands on its own two feet.

I will provide just such a theory in this section. Previous work makes this easier than it might at first seem. Elsewhere, I have shown how to axiomatize the Set Story,\footnote{See \textcite{Button:LT1}.} and the same techniques can be used for the Function Story, with only minor adjustments. (For readability, I relegate most of the fiddlier details to the appendices.)

\subsection{Notational preliminaries}
I begin with some notation. Where $\sigma$ and $\tau$ are any object-language terms, I write:
\begin{listbullet}
	\item``$\sigma(\tau)$'' to indicate the value of applying $\sigma$ to $\tau$.
	\item``$\isdefined{\tau}$'' to indicate that $\tau$ exists, i.e.\ that $\exists y\ y = 
	\tau$. 
	\item``$\sigma \simeq \tau$'' to indicate that $\sigma = \tau$ if either exists.\footnote{i.e.\ $\sigma \equalifdef \tau$ means $((\isdefined{\sigma} \lor \isdefined{\tau}) \lonlyif \sigma = \tau)$; here, $\sigma$ and $\tau$ are any object-language terms.}
	\item``$\Lamabs{x}{\tau}$'' for the (first-order) function which sends $x$ to $\tau$, if there is such a function.\footnote{i.e.\ $\Lamabs{x}{\tau}(v) \simeq \suballforin{\tau}{v}{x}\text{, for all }v$.}
\end{listbullet}
There are a few different ways to regiment these bits of notation more pedantically, but the exact choice of pedantry is relatively unimportant, so I relegate to \S\ref{s:pedantry:partial}. I will continue with a key definition, which allows me to relate functions to ``sets'':
\begin{define}\label{def:funset}
	We introduce abbreviations as follows:\footnote{As usual, $(\forall x \lhd y)\phi$ abbreviates $\forall x(x \lhd y \lonlyif \phi)$, for any infix predicate $\lhd$.} 
	\begin{align*}
		g \funin f &\colonequiv (\isdefined{f(g)} \lor  \exists y\ f(y) = g)\\
		g \funsubeq f &\colonequiv (\forall x \funin g)x \funin f
	\end{align*}
	Say that $f$ is a \emph{\funset} iff $(\forall x \funin f)f(x) = x$, i.e.\ $f$ is any partial identity function. \defendhere
\end{define}\noindent
We might read ``$g \funin f$'' as ``$g$ is in $f$'s field'' (i.e.\ in $f$'s domain or its range). We might read ``$g \funsubeq f$'' as ``$f$'s field includes $g$'s''. And \funsets will serve as our function-theoretic surrogates for sets, easily allowing us to construct a set-hierarchy as an inner model of a function-hierarchy.\footnote{Specifically, as in \S\ref{s:de:sketch}, we will interpret \ZF-sets as \emph{hereditary-\funsets}. But note that that interpretation does not witnesses the \emph{synonymy} between \ZF and \FLTzf (see the Main Theorem in \S\ref{s:de:sketch}). Under synonymy, \emph{every} \FLTzf-function ``simulates'' some \ZF-set, and \emph{every} \ZF-set ``simulates'' some \FLTzf-function (see footnote \ref{fn:IJonlybi}).
	
	Many alternative definitions of \funsets would have worked equally well. With \textcite[31]{Robinson:TCMVNS}, we could have defined \funsets as functions obeying the principle that $\forall x(\isdefined{f(x)} \lonlyif f(x) = 0)$. Or, with von Neumann (see footnote \ref{fn:vNimplementation}) we could have defined \funsets as characteristic functions. Nothing much turns on this choice (although the idea of a characteristic function requires that there be at least two objects, which neither \FLT or \FST guarantee). However, defining \funsets as partial identity functions gives an immediate link to category theory (see Theorem \ref{thm:flt:category} and \S\ref{app:FLT:cat}).}

\subsection{Using stages: \FST}
I now turn in earnest to the task of axiomatizing the Function Story. Since that Story mentions both stages and functions, I will start with a formal theory which does the same. My theory, then, has two distinct sorts of first-order variables, for \emph{functions} (lower-case italics) and for \textbf{stages} (lower-case bold). It has these two primitives:
\begin{listbullet}
	\item[$<$:] a relation between stages; read ``$\stage{r} < \stage{s}$'' as ``$\stage{r}$ is before $\stage{s}$''.
	\item[$\foundat$:] a relation between a function and a stage; read ``$f \foundat \stage{s}$'' as ``$f$ is found at $\stage{s}$''.
\end{listbullet}
For brevity, let $f \foundbefore \stage{s}$ abbreviate $\exists \stage{r}(f \foundat \stage{r} < \stage{s})$, i.e.\ $f$ is found before $\stage{s}$. Now let \FST, for Functional Stage Theory, be the theory with these five axioms:
\begin{listaxiom}
	\labitem{FunExt}{flt:ext} $\forall f \forall g\big(\forall x\, f(x) \equalifdef g(x) \lonlyif f = g\big)$	
	\labitem{FunOrd}{fst:ord} 
	$\forall \stage{r} \forall \stage{s}\forall \stage{t}(\stage{r} < \stage{s} < \stage{t} \lonlyif \stage{r} < \stage{t})$	
	\labitem{FunStage}{fst:stage} 
	$\forall f \exists \stage{s}\ f \foundat \stage{s}$
	\labitem{FunPri}{fst:pri} 
	$\forall \stage{s} (\forall f \foundat \stage{s})(\forall x \funin f) x \foundbefore \stage{s}$
	\labitem{FunSpec}{fst:spec} 
	$\forall P \forall \stage{s}\big(\forall x\big(\isdefined{P(x)} \lonlyif (x \foundbefore \stage{s} \land P(x) \foundbefore \stage{s})\big) \lonlyif \Lamabs{x}{P(x)}\foundat \stage{s}\big)$\\
	\emph{with ``$P$'' a second-order function-variable}
\end{listaxiom}
The first two axioms are analytic.  
\ref{flt:ext} says: functions are individuated by their mapping behaviour; I take it that this is a stipulative feature of the kind of functions we are considering.  
\ref{fst:ord} says: \emph{before} is a transitive relation on stage; I take it that this is required by the very idea of a notion of \emph{before}.\footnote{\label{fn:stagewo}Strikingly, there is no need to assume that \emph{before} is a well-order. Our axioms allow us to prove well-foundedness. We cannot prove linearity, but this is unimportant, since the ordinal height of a stage determines exactly which functions are found at that stage. This all follows from Theorems \ref{thm:fst-flt}--\ref{thm:flt:wo}, and it is characteristic of stage theories (cf.\ \cites[\S5]{Button:LT1}[\S4]{Button:LT2}[\S4]{Button:LT3}).}
	Then the next three axioms are clearly given by the Story. So: 
\ref{fst:stage} says: every function is found at some stage. 
\ref{fst:pri} says: everything in a function's field is found (strictly) before the function itself is found. And
\ref{fst:spec} says: if $P$ is a (second-order) map and everything in $P$'s field is found before $\stage{s}$, then $P$ determines a (first-order) function which is found at stage $\stage{s}$. 

I just invoked second-order logic. I will do so throughout, assuming (impredicative) second-order Comprehension. That said, my use of second-order logic is mostly just for notational convenience; almost everything I say can easily be first-orderized.\footnote{Specifically, \ref{fst:spec} would become this scheme. If $\forall x \forall y\forall z\big((\phi(x,y) \land \phi(x,z)) \lonlyif (y = z \land x \foundat \stage{s} \land y \foundat\stage{s})\big)$, then $(\exists f \foundat \stage{s})\forall x\forall y(f(x) = y \liff \phi(x,y))$.} The only case where second-order logic is required is in the discussion of quasi-categoricity (see \S\ref{s:quasi-categorical}).

Clearly, \FST's axioms are all true of the Function Story, as told in \S\ref{s:story}. So this is an excellent first step in axiomatizing that Story. But it is just a first step. Recall my aims from the start of \S\ref{s:flt}: I want to show that set-theoretic and function-theoretic approaches are on a par when it comes to mathematical foundations. Since we can axiomatize set theory without using a primitive \textbf{stage}-sort, I must next eliminate the \textbf{stage}-sort from function theory. 

\subsection{Eliminating stages: \FLT}\label{s:flt:formal}
With this goal in mind, let \FLT, for Functional Level Theory, be the theory with these three axioms:
\begin{listaxiom}
	\labitem{\ref{flt:ext}}{flt:ext2} $\forall f \forall g\big(\forall x\, f(x) \equalifdef g(x) \lonlyif f = g\big)$	
	\labitem{FunStrat}{flt:strat} $\forall a \exists s \big(\fevpred(s) \land a \funsubeq s\big)$
	\labitem{FunComp}{flt:comp} 
		$\forall P \big(\exists f \forall x\big(\isdefined{P(x)} \lonlyif (x \funin f \land P(x) \funin f)\big) \lonlyif \Lamabsalt{x}{P(x)}\text{ exists}\big)$
		\\\emph{with ``$P$'' a second-order function-variable}
\end{listaxiom}\noindent
As promised, these axioms do away with stages. But they require some discussion. We have already encountered \ref{flt:ext} in formulating \FST. Then \ref{flt:comp} is a principle which turns (second-order) maps into (first-order) functions; in effect, \ref{flt:comp} plays the same role in \FLT as \ref{fst:spec} played in \FST. The trickiest axiom to explain is \ref{flt:strat}. This uses a new predicate, ``$\fevpred$'', and we should pronounce ``$\fevpred(s)$'' as ``$s$ is a function-level'', or ``$s$ is a fevel'' for short. Intuitively, \FLT's fevels serve as function-theoretic proxies for \FST's stages, so that \ref{flt:strat} plays roughly the same role in \FLT as \ref{fst:stage} played in \FST. But, to ensure that the fevels behave this way, we must not take ``$\fevpred$'' as a primitive, but must instead define it. 

\emph{This can be done}. But the required definition is both cunning and ugly, and no one will be able to tell at a glance that it is the right definition. So I will not pause here to talk it through (see Definition \ref{def:flt:fevel} for more). What matters here is just that \emph{the definition works}.

The sense in which the definition \emph{works} is made precise by considering three characteristic theorems for \FLT. (For technical discussion, see \S\ref{app:FLT}; after announcing a Theorem, I list in brackets the axioms used in its proof):
	\begin{thm}\label{thm:fst-flt}
		$\FST \proves \phi$ iff $\FLT \proves \phi$, for any \FLT-sentence $\phi$.
	\end{thm}
	\begin{thm}[\ref{flt:ext} + \ref{flt:comp}]\label{thm:flt:wo}
		The fevels are well-ordered by $\funin$.\footnote{That is: (1) if some fevel is $F$, then there is an $\funin$-minimal fevel which is $F$; and (2) $s \funin t \lor s = t   \lor t \funin s$, for any fevels $s$ and $t$.}
	\end{thm}
	\begin{thm}[\FLT]\label{thm:flt:category}
		We have a category: the arrows are the functions; the objects and their identity arrows are exactly the \funsets; composition is functional composition.
	\end{thm}\noindent
Theorem \ref{thm:fst-flt} says that \FLT is equivalent to \FST, in the sense that they make exactly the same claims about functions. Theorem \ref{thm:flt:wo} says that our stage-surrogates, the fevels, are well-ordered by $\funin$. We can prove this result without using \ref{flt:strat}; when we add \ref{flt:strat}, we obtain the extremely powerful tool of $\funin$-induction (see Definition \ref{def:fevof}). Finally, Theorem \ref{thm:flt:category} says that our functions constitute a category (see \S\ref{app:FLT:cat}).

\subsection{\FLT is quasi-categorical}\label{s:quasi-categorical}
The Function Story introduces the (bare idea of the) iterative notion of \emph{function}. So, every cumulative hierarchy of functions answers to that Story. Since \FST's axioms are all true of that Story, every cumulative hierarchy of functions obeys \FST, and hence also \FLT, by Theorem \ref{thm:fst-flt}. Hence, by Theorems \ref{thm:flt:wo} and \ref{thm:flt:category}, 
\emph{every cumulative hierarchy of functions is a category which is arranged into well-ordered fevels}. This is a rather pleasant result.\footnote{Not least because the bare idea of the iterative notion of \emph{function}---as explicated via the Function Story---made no explicit reference to either categories or well-ordering (cf.\  \cites[210]{Scott:AST}[\S5]{Button:LT1}[\S4]{Button:LT2}[\S4]{Button:LT3}; also footnote \ref{fn:stagewo}).}

We can push this point further. The Function Story said nothing about the height of the hierarchy of functions. So, by design, \FLT says nothing about that either. Still, as a second-order theory, \FLT is quasi-categorical. Informally, we can spell out the point as follows: \emph{Any two cumulative hierarchies of functions are structurally identical for so far as they both run, but one may be taller than the other.}

There are at least two ways to explicate this idea, but \FLT is quasi-categorical on both explications: it is both \emph{externally} and \emph{internally} categorical. I will not go into details on this point here, since I have done so in detail elsewhere.\footnote{See \textcites[\S6]{Button:LT1}[\S5]{Button:LT3}.} However, an immediate upshot is that the size of a standard model of \FLT is determined entirely by the order-type of its fevels. Specifically, consider this function, $p$, from ordinals to cardinals:
\begin{align*}
	p(1) &\coloneq 1 & 
	p(\alpha + 1)&\coloneq (p(\alpha)+1)^{p(\alpha)} & 
	p(\alpha)&\coloneq \sup_{\beta < \alpha}p(\beta)\text{ for limit }\alpha
\end{align*}
For any ordinal $\alpha> 0$, there are (standard) models of \FLT with (exactly) an $\alpha$-sequence of fevels. And, given an $\alpha$-sequence of fevels, there are $p(\alpha)$ functions.\footnote{\emph{Proof sketch.} By induction; the only interesting clause is $p(\alpha+1)$. Here, note that there are exactly $|B|^{|A|}$ \emph{total} functions $A \functionto B$; so, allowing for the one extra possibility of being undefined rather than landing in $B$, there are exactly $(|B|+1)^{|A|}$ \emph{partial} functions $A \functionto B$.} Evidently $p$ grows quite quickly: one fevel yields one function;\footnote{Specifically, the null function, which is undefined for all inputs.} but with four fevels we already have a billion functions.

\subsection{Adding more functions}
\FLT provides us with a mathematical playground of functions, arranged in well-ordered fevels. However, since \FLT allows that there can be only one fevel, and hence only one function, the vanilla theory of \FLT allows our playground to be pretty desolate. 

To ensure that our hierarchy of functions has a certain height---to populate the playground---here are some principles that we could add to \FLT:
\begin{listaxiom}
	\labitem{FunEndless}{flt:cre} $\forall s \exists t\ s \funin t$
	\labitem{FunInfinity}{flt:inf} $\exists s(\exists p \funin s)(\forall q \funin s)\exists r(q \funin r \funin s)$
	\labitem{FunSupercomp}{flt:rep} $\forall P\big(\exists f \forall x\big(\isdefined{P(x)} \lonlyif x \funin f\big) \lonlyif \Lamabsalt{x}{P(x)}\text{ exists}\big)$
\end{listaxiom}\noindent
Against the background of \FLT, we can characterize these principles as follows. \ref{flt:cre} says: every fevel has a successor (compare \ZF's Powersets). \ref{flt:inf} says: there is an $\omega^\text{th}$ fevel (compare \ZF's Infinity). And \ref{flt:rep} strictly strengthens \ref{flt:comp}, saying: no function's field can be mapped unboundedly into the fevels (compare \ZF's Replacement).

\subsection{Synonymy}\label{s:de:sketch}
Let \FLTzf be the theory \FLT + \ref{flt:inf} + \ref{flt:rep}. The following result explains the name:\footnote{I have formulated \FLTzf as a \emph{second-order} theory; then the result is that \FLTzf is synonymous with \emph{second-order} \ZF. However, we also have synonymy between the \emph{first-order} versions of both theories.}
\theoremstyle{generic}
\newtheorem*{mainthm*}{\normalfont\bfseries{Main Theorem}}
\begin{mainthm*}
	\ZF and \FLTzf are synonymous.\footnote{\label{fn:Ipreferltplus}In fact, we can prove a much stronger result. Let \FLTplus be the theory \FLT + \ref{flt:cre}, i.e.\ the theory which arises by adding ``there is no last stage'' to the Function Story. Comparably, \LTplus is the theory which arises by adding ``there is no last stage'' to the Set Story (see \S\ref{s:app:lt} for details; notably, \LTplus is a strict sub-theory of \ZF). Then the stronger version of the Main Theorem is this: \emph{\LTplus and \FLTplus are synonymous}. 
		
	Since \ZF is the ``industry standard'' of set theory, I will only discuss the case of \ZF and \FLTzf in the main text. But I would much prefer to focus on the case of \FLTplus; after all, focussing on the weaker theories only makes my argument stronger. (I revisit this in footnotes \ref{fn:coherentltplus}, \ref{fn:ohweaker}, and \ref{fn:secondplus}.)} 
\end{mainthm*}\noindent
I outlined the notion of \emph{synonymy} in \S\ref{s:vonNeumann}. But here is a helpful, if rough, gloss of the Main Theorem: in \emph{some} sense, \ZF and \FLTzf are notational variants; they wrap the same content in different notational packaging. 

This rough gloss is deliberately provocative, and it raises many philosophical questions. I will explore those questions in the next two sections. First, I should quickly sketch the proof strategy for the Main Theorem.

In \ZF, we can regard functions as \setfunctions, i.e.\ sets of ordered pairs of a certain sort (see \S\ref{s:vonNeumann}). This provides us with an interpretation $I : \FLTzf \functionto \ZF$. In particular: working in \ZF, we can interpret \FLTzf as dealing with \emph{hereditary-\setfunctions}, where a set is a hereditary-\setfunction iff it is a \setfunction, whose field comprises only \setfunctions, whose field comprises only \setfunctions, etc.\footnote{\label{fn:IisInt}Showing that this \emph{is} an interpretation takes more effort than one might expect (particularly if the aim is to establish the results for \LTplus and \FLTplus; see footnote \ref{fn:Ipreferltplus}). The key steps are as follows: first, prove \ref{flt:ext}$^I$ and \ref{flt:comp}$^I$ within \LTplus; then you can appeal to the $I$-translation of (all the results you used to establish) Theorem \ref{thm:flt:wo}; these entail the good behaviour of the fevels$^I$, thereby delivering \ref{flt:strat}$^I$. (Similar but simpler moves are made in \cite[\S{}C]{Button:LT3}.)}

Conversely, in \FLTzf, we can interpret sets as \funsets (see Definition \ref{def:funset}). This provides us with an interpretation $J : \ZF \functionto \FLTzf$. In particular: working in \FLTzf, we can interpret \ZF as dealing with \emph{hereditary-\funsets}, where a function is a hereditary-\funset iff it is a \funset, whose members (in the sense of $\funin$) are all \funsets, whose members are all \funsets, etc.\footnote{Just as we can use transfinite recursion in \ZF to define \emph{hereditary-\setfunctions}, we can use transfinite recursion in \FLTzf to define \emph{hereditary-\funsets}. The proof that $J$ \emph{is} an interpretation is similar to the proof that $I$ is an interpretation (see footnote \ref{fn:IisInt}).}

We now have interpretations going in both directions. With a little effort,\footnote{The key tools are: patience, $\in$-induction in \LTplus, and $\funin$-induction in \FLTplus.} we can show that they are \emph{bi}-interpretations. By the Friedman--Visser Theorem,\footnote{\label{fn:IJonlybi}i.e.\  \textcite[Corollary 5.5]{FriedmanVisser:WBIS}. Note that the interpretations $I$ and $J$ do not \emph{themselves} bear witness to the synonymy. Rather, the Friedman--Visser Theorem shows that we can use $I$ and $J$ to define other interpretations which witness synonymy.} it follows that \ZF and \FLTzf are synonymous, as required.

\section{Function-theoretic judicial foundations}\label{s:judicial}
Let me take stock. I began \S\ref{s:story} with the Function Story, which introduced the bare idea of the iterative notion of \emph{function}. The Story is quasi-categorically axiomatized by \FLT; so the Story entails that functions constitute a category which is arranged into well-ordered fevels. Moreover, we can easily enrich \FLT with various claims to ensure that the sequence of fevels runs for a rather long time. In particular, the theory \FLTzf arises just by insisting that the hierarchy of functions is extremely tall. And \FLTzf is synonymous with \ZF. 

This synonymy suggests that we might be able to replace set-theoretic foundations (as formalized in \ZF) with function-theoretic foundations (as formalized in \FLTzf). In this section, I will show that this is correct. Indeed, I will show that function-theoretic and set-theoretic foundations provide us with exactly the same (judicial) foundations.

\subsection{Maddy on judicial foundations}\label{s:judicialfoundations}
Many things have been meant by talk of mathematical ``foundations''. So I must start by clarifying the kind of foundational project I have in mind. In brief: my interest is not in \emph{methodological} foundations, nor (initially) in \emph{metaphysical} foundations, but in (what I call) \emph{judicial} foundations. 

\
\\
\textbf{Methodological foundations.}
A mathematical theory might aim to supply us \emph{methodological} foundations. Specifically, it might aim to provide a vocabulary and toolkit for dealing with situations that arise frequently across the motley of mathematics. 

Category theory has some claim to providing methodological foundations; think of the unifying power of concepts like \emph{functor}, \emph{natural transformation}, or \emph{adjunction}. Set theory, by contrast, has almost no such claim: to give a well-worn example, having shown that we can define the real numbers as sets of sets of sets of\ldots sets, no practising mathematician ever digs back through that definition. Function theory has a similarly weak claim to providing methodological foundations, and for similar reasons. So I will not discuss methodological foundations any further. 

\
\\
\textbf{Metaphysical foundations.}
Another possible foundational project is to determine what mathematical objects \emph{are}. 

Those who ask such questions are not usually satisfied with a grab-bag answer, along the lines of: \emph{mathematical objects include numbers, groups, graphs, topologies, categories,\ldots}.  
They will reply: \emph{yes, but what exactly \textsc{are} these things?} And they will not be satisfied with a mathematical definition, like: \emph{a category is a system of objects and arrows such that\ldots}. They are in pursuit of something distinctively metaphysical. 

Such people are asking for a \emph{metaphysical foundation} for mathematics. I will set this metaphysical project aside for now, returning to it in \S\ref{s:metaphysics}, where we will be able to shed new light on it.  

\
\\
\textbf{Judicial foundations.} My focus in this section, instead, is on \emph{judicial foundations}. The inspiration behind this phrase is in Penelope Maddy's suggestion that set theory is ``a court of final appeal''.\footnote{\textcite[26]{Maddy:NM}. \textcite[296]{Maddy:STF} has since expressed worries that the notion of ``a final court of appeal'' was ``something of an exaggeration'', and has outlined various nuanced notions in that ballpark. I would particularly emphasize the continuity between what I say about ``judicial foundations'' and \citepossess[289--98]{Maddy:STF} notions of a Meta-mathematical Corral, a Generous Arena, and a Shared Standard.} Here is her point, in detail: 
\begin{quote}
	The force of set theoretic foundations is to bring (surrogates for) all mathematical objects and (instantiations of) all mathematical structures into one arena---the universe of sets---which allows the relations and interactions between them to be clearly displayed and investigated. Furthermore, the set theoretic axioms developed in this process are so broad and fundamental that they do more than reproduce the existing mathematics; they have strong consequences for existing fields and produce a mathematical theory that is immensely fruitful in its own right. Finally, perhaps most fundamentally, this single, unified arena for mathematics provides a court of final appeal for questions of mathematical existence and proof: if you want to know if there is a mathematical object of a certain sort, you ask (ultimately) if there is a set theoretic surrogate of that sort; if you want to know if a given statement is provable or disprovable, you mean (ultimately), from the axioms of the theory of sets.\footnote{\textcite[26]{Maddy:NM}.}
\end{quote}
The crucial point is something like this.\footnote{In addition to taking this from Maddy's writings, I would like to thank Michael Potter for impressing this point upon me many years ago.} Suppose you make a conjecture, in any area of mathematics. If someone can construct a set, in \ZF, which provides a counter-example to the conjecture, then you lose. If you can show that the conjecture holds in (vanilla) \ZF, then you win. And if you can show that the conjecture holds in \ZF plus some extra large cardinal axiom, people might well be interested.

This role for \ZF is not, I think, a mere sociological quirk concerning mathematical practice. Instead, I think there are two very good reasons for regarding standard \ZF as an excellent candidate for judicial foundations: 
\begin{listl-0}
	\item\label{n:things} it provides us with plenty of objects, which are richly structurable; and
	\item\label{n:consistent} it is consistent and indeed coherent.
\end{listl-0}
I will explain what I mean by ``coherent'' in \S\ref{s:foundation:consistent}. But let me quickly unpack these points, starting with \eqref{n:things}. \ZF tells us that there are plenty of objects. Moreover, the membership relation arranges these objects in an extremely tractable way. Using this given membership-structure, we have the resources to define many further relations on our objects, structuring and restructuring them to our heart's content. This is the point of \eqref{n:things}. Indeed, this is \emph{why} (in Maddy's phrase) we have a rich paradise of ``set-theoretic surrogates'' for other mathematical objects. Moreover, set theory gives us these objects, with their rich structure, without exploding under its own weight. That is the point of \eqref{n:consistent}.

Features \eqref{n:things} and \eqref{n:consistent}, then, are our two hallmarks for a decent judicial foundation. I will further elaborate on both features below, but my basic task for this section is to show that function theory possesses both features, to exactly the same extent as set theory. (Note: in \S\S\ref{s:foundation:consistent}--\ref{s:foundationsketch}, I will not really distinguish ``set theory'' from \ZF, or ``function theory'' from \FLTzf; I will explain why this is acceptable \S\ref{s:tiedforfuture}.) 

\subsection{Consistency and coherence}\label{s:foundation:consistent}
I begin by considering \eqref{n:consistent}. I will argue that both set theory and function theory are consistent and coherent, starting with set theory. 

Here is the strongest argument that I know of for the consistency of \ZF:
\begin{listbullet}
	\item[(1$_\in$)]\label{arg:bare} The Set Story of \S\ref{s:story} introduces the bare idea of the iterative notion of \emph{set}. We can augment that Story with two further claims: there is some limit stage; and no set can be mapped unboundedly into the sequence of stages.	Call the result the Augmented Set Story. 
	\item[(2$_\in$)] The Augmented Set Story is consistent. After all, the (unaugmented) Set Story is obviously consistent, and\footnote{\label{fn:coherentltplus}This is obviously the weakest point of this argument. How can we be sure that adding these two claims doesn't undermine consistency? If this really worries us, we can retreat from \ZF to \LTplus (see footnote \ref{fn:Ipreferltplus} and \S\ref{s:app:lt}). The theory \LTplus quasi-categorically axiomatizes the Set Story when it is augmented only with the claim that there is no last stage \parencite[see][]{Button:LT1}, and I literally cannot entertain the idea that such a weakly-augmented story is inconsistent. Comparably, \FLTplus arises by augmenting the Function Story (only) with the claim that there is no last stage. Moreover, \LTplus and \FLTplus are synonymous. So the entire argument of this subsection can be rerun in terms of \LTplus and \FLTplus rather than \ZF and \FLTzf, and the argument will be all the stronger for it.} the two further claims are also consistent.
	\item[(3$_\in$)]\label{arg:axiom} \ZF quasi-categorically axiomatizes the Augmented Set Story.\footnote{For this, see \textcite{Button:LT1}.} So \ZF inherits that Story's consistency.
\end{listbullet} 
I think this argument is a good one; indeed, it is why I believe that \ZF is consistent. Moreover, the argument generalizes in an important way. 

To explain the generalization, note that consistency is a rather weak property. Let $\Thnamed{PA}_*$ be $\Thnamed{PA} + \lnot\textspaced{Con}(\Thnamed{PA})$.\footnote{Here, ``$\textspaced{Con}(\ldots)$'' indicates a canonical way for forming consistency sentences; I am assuming(!)\ \Thnamed{PA}'s consistency.} By G\"odel's second incompleteness theorem, $\Thnamed{PA}_*$ is consistent. But $\Thnamed{PA}_*$ cannot be ``trusted''; after all, $\Thnamed{PA}_*$ ``proves its own inconsistency'', in the sense that $\Thnamed{PA}_* \proves \lnot\textspaced{Con}(\Thnamed{PA}_*)$. So, what we want is not just that our theory should be \emph{syntactically} consistent, but that---speaking informally---our theory should have an \emph{intended interpretation}.\footnote{\label{fn:coherentinterpretation}I do not want to assume that the notion of ``interpretation'' must be \emph{model-theoretic}, for two reasons. First, I would then have to decide whether to formalize model theory set-theoretically or function-theoretically. Second, for reasons I have described at length elsewhere \parencites[chs.\ 6--12]{ButtonWalsh:PMT}{Button:MIR}, I would rather work within a second-order object language than within a first-order metalanguage.} Following Stewart Shapiro, I refer to this as the demand that the theory should be \emph{coherent}.\footnote{\textcite[95--6, 133--6]{Shapiro:PM}.}

Fortunately, and as hinted, the argument (1$_\in$)--(3$_\in$) can be strengthened. The Set Story is not presented as uninterpreted syntax; it is presented in English, with an associated intended interpretation. Indeed, we are convinced that the Story is consistent precisely because it strikes us as coherent.\footnote{\label{fn:coherentisaacson}The obvious question arises: why should we think there are (or could be) enough objects to provide us with this intended interpretation of \ZF? Read a certain way, this question will swiftly lead us from judicial into metaphysical foundations, whereupon I do not know exactly how to respond (see \S\ref{s:metaphysicsnotdead}). However, I do want to signal my agreement with \textcite[29]{Isaacson:RMCST}, that we take a big step towards showing that \ZF is coherent when we show that \ZF quasi-categorically axiomatizes the Augmented Story, thereby ``developing our mathematical understanding of the subject matter of this theory through informal rigor''.} So we can rerun the argument (1$_\in$)--(3$_\in$), with ``coherent'' in place of ``consistent''. 

This concludes my case that set theory is both consistent and coherent. With this to hand, I now want to show the same of \FLTzf. 

Here is a cute argument. Since \ZF and \FLTzf are synonymous, they are equiconsistent. We just saw that \ZF is consistent. So \FLTzf is consistent too. QED? 

The cute argument has two weaknesses. First: it is not completely obvious that synonymy always preserves coherence (rather than mere consistency). Second: if this is our \emph{only} argument for \FLTzf's consistency, then \FLTzf's consistency is epistemically dependent upon \ZF's; in which case, at least epistemically, set theory rather than function theory will serve as the ``ultimate'' court of appeal. 

Fortunately, both issues can be dealt with. An almost identical argument, with exactly the same suasive force, establishes the consistency of \FLTzf: 
\begin{listbullet}
	\item[(1$_@$)] The Function Story introduces the bare idea of the iterative notion of \emph{function}. We can augment that Story with two further claims: there is some limit stage; and no function's field can be mapped unboundedly into the sequence of stages. Call the result the Augmented Function Story. 
	\item[(2$_@$)] The Augmented Function Story is consistent. After all, the (unaugmented) Function Story is obviously consistent, and the two further claims are consistent.
	\item[(3$_@$)] \FLTzf quasi-categorically axiomatizes the Augmented Function Story. So \FLTzf inherits that Story's consistency.
\end{listbullet} 
Moreover, as before, we can rerun the argument with ``coherent'' in place of ``consistent''.\footnote{With, mutatis mutandis, the clarifications of footnotes \ref{fn:coherentltplus}, \ref{fn:coherentinterpretation} and \ref{fn:coherentisaacson}. Indeed, this highlights part of the importance of \FLT's quasi-categoricity (see \S\ref{s:quasi-categorical}).} So \FLTzf is both consistent and coherent.

\subsection{Function-theoretic foundations}\label{s:foundationsketch}
I now turn to \eqref{n:things}. In \S\ref{s:judicialfoundations}, I explained the sense in which \ZF provides us with a vast array of richly structurable objects. I will now show that the same holds of \FLTzf. 

Once again, synonymy supplies us with a cute argument to this effect. Specifically: there is (ultimately) a \ZF-surrogate for an object iff there is (ultimately) a \FLTzf-surrogate for an object; just apply the synonymy to the description of the surrogate. QED.

Again, one might worry that this cute argument threatens to diminish function theory's status. After all: if our \emph{only} method for finding surrogate-objects involves first using set theory and then invoking synonymy, then set theory may seem to have a kind of primacy over function theory.  

This time, though, the worry is wholly misguided: it confuses \emph{judicial} foundations with \emph{methodological} foundations. Granted, in searching for a surrogate-object, sometimes set theory will be more wieldy than function theory. Moreover, issues of (un)wieldiness are of real practical importance. But they are orthogonal to judicial foundations. In judicial foundations, the question is whether there are \emph{ultimately} surrogates for objects or derivations, and the phrase ``ultimately'' is allowed to conceal all manner of unwieldiness. Indeed, from a purely judicial perspective, it would not matter if we presented function-theoretic foundations just by starting with our favourite set-theoretic foundations and then translating them into function theory using the interpretation $J : \ZF \functionto \FLTzf$ from \S\ref{s:de:sketch}. The result will be less wieldy than set-theoretic foundations we started with; still, we can (ultimately) find function-theoretic surrogates iff we can (ultimately) find set-theoretic surrogates, and that is the only thing which matters, judicially speaking. 

(All that said, I should probably note that there is no need to make function-theoretic foundations ``parasitic'' upon set-theoretic foundations. We can present ``autonomous'' function-theoretic foundations, providing surrogates for various mathematical entities in a setting where our basic unit of currency is functions rather than sets. Especially enthusiastic readers may enjoy playing around with this possibility for themselves, and they will can find some playful ideas in \S\ref{app:Judicious}. But my philosophical ambitions do not require me to force less enthusiastic readers to dwell on this point.)

\subsection{Adding new axioms}\label{s:tiedforfuture}
So far, I have argued that \FLTzf meets conditions \eqref{n:things} and \eqref{n:consistent} exactly as well as \ZF. 
But I have spoken as if \emph{function theory} meets \eqref{n:things} and \eqref{n:consistent} exactly as well as \emph{set theory}. It may not be immediately obvious that the conflation between \ZF and set theory is harmless. After all: set theorists often consider extensions of \ZF,\footnote{\label{fn:ohweaker}Of course, set theorists also sometimes consider theories which are \emph{weaker} than \ZF. As per footnotes \ref{fn:Ipreferltplus} and  \ref{fn:coherentltplus}, I can re-run the arguments of the last few sections using just \LTplus.} and presumably function theorists (if there were any) would often consider extensions of \FLTzf. 
%
So, for all I have said so far, one might worry that the following situation is possible: 
set-theoretic considerations push us to add some new axiom, $\phi$, to \ZF; 
function-theoretic considerations push us to add some new axiom, $\psi$, to \FLTzf; 
but $\ZF + \phi$ is \emph{not} synonymous with $\FLTzf + \psi$, and indeed our favoured set-theoretic foundations and our favoured function-theoretic foundations have come apart.

Roughly put, my reply to this worry is as follows. Considering synonymy will show that we never really have a case where set-theoretic and function-theoretic considerations pull us in different directions in the search for new axioms. To refine this reply, though, I will separately consider extrinsic and intrinsic arguments for new axioms.\footnote{See e.g.\ \textcites{Maddy:BA1}{Maddy:BA2}.}

\
\\
\textbf{Extrinsic arguments for new axioms.} An extrinsic argument for $\phi$ is an argument that we should accept $\phi$ because doing so has better all-things-considered consequences for other areas of mathematics. I claim that synonymy always preserves extrinsic arguments. 

To see why, suppose that the synonymy is witnessed by translations $\intTtoS: \FLTzf \functionto \ZF$ and $\intStoT: \ZF \functionto \FLTzf$.\footnote{\label{fn:secondplus}As in footnote \ref{fn:Ipreferltplus}, this argument can be run using \LTplus and \FLTplus.} Then $\ZF + \phi$ is synonymous with $\FLTzf + \phi^\intStoT$, for every set-theoretic sentence $\phi$. Similarly, $\ZF + \psi^\intTtoS$ is synonymous with $\FLTzf + \psi$, for every function-theoretic sentence $\psi$. Now, synonymous theories (ultimately) have the same consequences for all areas of mathematics. So, using ``best'' as short-hand for ``(all things considered) ultimately best for mathematics'', we have two schematic biconditionals:
\begin{listbullet}
	\item[(1E)] It is best to add $\phi$ to \ZF iff it is best to add $\phi^\intStoT$ to \FLTzf. 
	\item[(2E)] It is best to add $\psi$ to \FLTzf iff it is best to add $\psi^\intTtoS$ to \ZF. 
\end{listbullet}
Now consider a situation where we are worried that set-theoretic and function-theoretic considerations are, for extrinsic reasons, pulling us in different directions. For example, and without loss of generality, suppose that Anne claims that it is best to add \CH to \ZF, and Bob claims that it is best to add (some statement which entails) $\lnot\CH^\intStoT$ to \FLTzf. By (1E), Bob must accept that it is best to add $\lnot\CH$ to \ZF. But then it is clear that the disagreement between Anne and Bob has nothing per se to do with sets \emph{versus} functions. Anne might have had essentially the same disagreement---about whether to accept or reject the continuum hypothesis---with a set theorist who has never heard of \FLTzf. And exactly the same consequences, for all areas of mathematics, would have been at stake in the purely intra-set-theoretic formulation of this disagreement.

\
\\
\textbf{Intrinsic arguments for new axioms.} 
Having considered extrinsic arguments, let me now consider intrinsic arguments for new axioms. An intrinsic argument for $\phi$ is an argument that we should accept $\phi$ because $\phi$ is something like a conceptual truth. Again, I will claim that such arguments are preserved under synonymy, in the specific case of intrinsic arguments concerning the iterative notions of \emph{set} and \emph{function}. 

To argue for this claim, I want to start by establishing these two conditionals:
\begin{listr-0}
	\item[(1i)] If \ZF intrinsically supports $\phi$, then \FLTzf intrinsically supports $\phi^J$.
	\item[(2i)] If \FLTzf intrinsically supports $\psi$, then \ZF intrinsically supports $\psi^I$.
\end{listr-0}
Here, $I: \FLTzf \functionto \ZF$ and $J : \ZF\functionto \FLTzf$ are the translations mentioned in \S\ref{s:de:sketch}. So $I$ renders \FLTzf-functions as hereditary-\setfunctions, and $J$ renders \ZF-sets as hereditary-\funsets.

To establish (1i), suppose we are working within \FLTzf. Suppose, too, we read all mention of ``sets'', in the Augmented Set Story of \S\ref{s:foundation:consistent}, as concerned with sets$^J$, i.e.\  hereditary-\funsets (as defined in \S\ref{s:de:sketch}). Then we will agree that $\ZF^J$ is a quasi-categorical axiomatization of the Augmented Set Story (read in this way). So, whatever reasoning is supposed to show that \ZF intrinsically supports $\phi$, that same reasoning will show (to an exactly equal extent) that $\ZF^J$ supports $\phi^J$. Finally, note that \FLTzf proves $\ZF^J$. So, the intrinsic argument for $\phi$ becomes an intrinsic argument for $\phi^J$, as required. An exactly similar argument establishes (2i).

Using the conditionals (1i) and (2i), I now want to turn them into biconditionals:
\begin{listr-0}
	\item[(1I)] \ZF intrinsically supports $\phi$ iff \FLTzf intrinsically supports $\phi^J$.
	\item[(2I)] \FLTzf intrinsically supports $\psi$ iff \FLTzf intrinsically supports $\psi^I$.
\end{listr-0}
The left-to-right direction of (1I) is just (1i). For the right-to-left direction, suppose \FLTzf intrinsically supports $\phi^J$. Then \ZF intrinsically supports $\phi^{JI}$, by (2i). Since $I$ and $J$ comprise a bi-interpretation, $\ZF \proves \phi \liff \phi^{JI}$. So \ZF intrinsically supports $\phi$, delivering (1I). An exactly similar argument establishes (2I). 

Armed with (1I) and (2I), we can now argue as we did concerning extrinsic arguments. For example, and without loss of generality, suppose that Anne claims that \ZF intrinsically supports \CH, and Bob claims that \FLTzf intrinsically supports $\lnot\CH^J$. By (1I), Bob must accept that \ZF intrinsically supports $\lnot\CH$. And so, again, Anne and Bob's disagreement has nothing per se to do with sets \emph{versus} functions. 

\
\\
To be clear: I am not denying that Anne may not find some heuristic value in thinking in terms of sets rather than functions (or vice versa, for Bob).\footnote{Indeed, the above suggests a minor variant on the approach of seeking intrinsic arguments for new axioms: tell the Function Story (or something similar); axiomatize it with a theory (like \FLTzf) which is synonymous with \ZF; arguing that (some natural extension of) the Story intrinsically motivates $\phi$; infer that $\phi$-under-interpretation holds of the sets.} And I have taken no stance over which of Anne and Bob is \emph{right}; I am not suggesting any method for \emph{deciding} who is right; I am not even offering an account of what it might even \emph{mean} to be right. My point is only that, by considering synonymy, we can see that Anne's and Bob's disagreement is not one where set-theoretic foundations rival function-theoretic foundations \emph{per se}. There may be genuine rivalry over \CH; but iff so, there is exactly the same rivalry over $\CH^\intStoT$.

\subsection{The same judicial foundation}
Allow me to recap the main points of this section. Set theory coherently presents us with an astonishingly fertile array of richly structurable objects; objects which can serve as surrogates for mathematical objects from almost anywhere else. This is the sense in which set theory provides us with a mathematical paradise. 

Exactly the same, though, can be said of function theory. Function theory equally coherently introduces us to an equally fertile array of richly structurable objects. It presents us with an alternative but equally splendid vision of paradise. 

Indeed, we can now appreciate that these are two visions of the \emph{very same} (judicial) paradise. Set theory and function theory establish the same ultimate court of final appeal---they supply the very same judicial foundation---since all their verdicts are equally trustworthy (given \S\ref{s:foundation:consistent}) and ultimately identical (by \S\S\ref{s:foundationsketch}--\ref{s:tiedforfuture}). The difference between set-theoretic and function-theoretic judicial foundations amounts, not to varying the court, but simply varying the language within which the judicial proceedings are conducted.

\section{Joint-carving metaphysical foundations}\label{s:metaphysics}
So say I. But some will doubt that set-theoretic and function-theoretic foundations are really \emph{identical}. They will insist come that set theory and function theory postulate \emph{rival ontologies}: a paradise of sets, on the one hand, or a paradise of functions, on the other. Addressing this idea will move us away from considering judicial foundations, and back to metaphysical foundations (which I set aside in \S\ref{s:judicialfoundations}).

In this section, I will examine and reject
\emph{joint-carving realism} about metaphysical foundations. (In \S\ref{s:metaphysicsnotdead}, I will consider whether there can be metaphysical foundations without joint-carving.) The notion of joint-carving that I have in mind is encapsulated by Ted Sider's call-to-arms:
\begin{quote}
	The world has a distinguished structure, a privileged description. For a representation to be fully successful, truth is not enough; the representation must also use the right concepts, so that its conceptual structure matches reality’s structure. There is an objectively correct way to ``write the book of the world''.\footnote{\textcite[i]{Sider:WBW}.}
\end{quote}
The case study of this paper raises an interesting challenge for joint-carving realists (a challenge which builds on Putnam; see \S\ref{s:Putnam} for more). In brief: since set theory and function theory are synonymous, both can be made true if either can. But, according to joint-carving realists, for a theory ``to be fully successful, truth is not enough''; a \emph{fully} successful theory must ``use the right concepts''. Joint-carving realists therefore face a substantial question: is set-\emph{membership} or functional-\emph{application} one of ``the right concepts''?\footnote{\textcite[\S7.13, \S8.2.2, \S8.3.3, ch.13n.2]{Sider:WBW} himself often suggests that membership might be joint-carving. Given the possibility of function-theoretic foundations, he should reconsider whether, instead, application is joint-carving.} 

In this section, I will argue that it is a disaster to think that this question has an answer. First, though, I should flag a general discomfort with joint-carving realism. Joint-carving realists ask whether certain concepts are fundamental, where concepts are the semantic analogues of \emph{predicates} and \emph{function-symbols}.\footnote{A few lines below his call-to-arms, \textcite[i]{Sider:WBW} writes: ``a fact is fundamental when it is stated in joint-carving terms\ldots. Questions about which expressions carve at the joints are questions about how much structure reality contains.''} A priori, it is strange to think that this question should be asked. Perhaps an idealist should expect that the structure of the (phenomenal) world must be grammatically tractable, so that ``the joints of nature'' would align with a specially chosen sub-sentential vocabulary. But I am not an idealist, and nor are joint-carving realists like Sider. I see no reason whatsoever to think that reality must have a fundamental structure which mirrors the grammatical structure of recognisably human languages. A fortiori, I feel no compulsion whatsoever to ask whether set-theoretic \emph{membership} or function-theoretic \emph{application} is joint-carving. 

Still, the question \emph{has} been posed. My aim is to argue that it should not be answered. 

\subsection{Monistic joint-carving and its surds}\label{s:Surd}
It will help to make things personal. Allow me to introduce two characters, setty-Seb and functiony-Fern. They are joint-carving realists who adopt rival positions on the question at hand. Here they are, in stereo:
\begin{quote}
	\begin{center}
		\begin{tabular}{@{}p{0.41\textwidth}p{0.41\textwidth}@{}}
			
			\textbf{Seb:} Membership is joint-carving.	
			Application is derivative, induced from membership by interpretation $I$. 
			&
			\textbf{Fern:} Application is joint-carving. Membership is derivative, induced from application by interpretation $J$.
			\tabularnewline\addlinespace
%
%
%
		\end{tabular}
	\end{center}
\end{quote}
Seb and Fern are offering \emph{rival} approaches, in the sense that their claims cannot consistently both be advanced in one and the same breath. The joint-carving realist must, then, hold that at most one of Seb and Fern is right. 

Now, one might think that things have already gone wrong, on the grounds that all synonymous theories must be equally joint-carving. But that would be mistaken. Not only is there no good general argument that synonymy preserves joint-carving-ness; synonymy does not always preserve (intuitive) truth. This is because, whilst synonymy establishes that two theories are neatly inter-translatable, the translations witnessing synonymy might be \emph{mis}translations. To illustrate: a child who systematically confuses the words ``pug'' and ``French bulldog'' may offer a theory which is synonymous with mine (just translate ``pug'' to ``French bulldog'' and vice versa); they still say something false when, pointing at my handsome hound, Ragnar, they say ``that's a pug''.\footnote{Of course, it is not always a given that children intend to speak the same language as those around them; they may intend to speak a code. The point is: if they intend to speak the same language as their peers, then (ceteris paribus) the correct translation is homophonic, so they say something false.}

That said, my firm insistence that Ragnar is a French bulldog (not a pug) has very little in common with Seb's insistence that membership (not application) is fundamental. For one thing, we can point at dogs, but not at sets or functions (let alone at the metaphysician's idea of \emph{joint-carving}). For another, as we saw in \S\ref{s:tiedforfuture}, there is no \emph{mathematical} disagreement between Seb and Fern. Their apparent disagreement is purely \emph{metaphysical}. And, as I will now argue, it is absurd. 

To see that Seb and Fern's disagreement is absurd, let us start by considering how it might be settled. Very likely, Seb, Fern, and any other joint-carving realists will point us in the direction of \emph{metaphysical virtues}. We will be told to evaluate matters on the basis of: explanatoriness; unificatoriness; simplicity; tractability; ideological/ontological parsimony; etc. 

For the sake of argument, I will make two huge concessions to the joint-carving realist: first, that these virtues are well-defined; second, that greater virtue is a reliable indicator of greater joint-carving-ness.\footnote{For the record: I suspect the first is false and the second is nonsense.} Even making these concessions, Seb's and Fern's approaches score equally well as regards any metaphysical virtues. 

To be sure, in different contexts, pragmatic considerations may favour one approach over the other. Perhaps the functions-first approach is more immediately hospitable to certain kinds of algebraic reasoning; perhaps the sets-first approach is more immediately amenable in other settings. But as we saw in \S\ref{s:judicial}, there is no real mathematical breathing-room between the two approaches. Consequently, it is hopeless to believe that metaphysical virtues will help to decide which of sets and functions is fundamental. 

Seb can, of course, simply stamp his feet and insist that membership simply \emph{is} more fundamental than application. (Fern can stamp her feet right back at him.) But he will have to admit that there is no \emph{reason why} this is so;\footnote{So Seb would be violating any version of the principle of sufficient reason; see \textcite{Amijee:PSR} for a discussion of its resurgence in modern metaphysics.} it is just a \emph{brute fact}, an inexplicable and unknowable \emph{metaphysical surd}.\footnote{I am deliberately using \citepossess[46--8]{Putnam:RTH} phrase; see also \S\ref{s:Putnam}.} Positing such a surd is formally consistent, of course, but it is extremely unappealing. Indeed, in the cost/benefit game that metaphysicians often play,\footnote{The idea is traceable to \textcite[14]{Quine:OI}; but \textcite[4]{Lewis:OPW} gave us the method for playing it.} such a surd is a substantial cost. 

The costs continue to rise, as we move from considering metaphysics to metasemantics. To see why, let us try to imagine (without loss of generality) that Seb is right to stamp his feet, because {membership} rather than {application} is joint-carving. So, according to Seb: \emph{there are a great many objects, fundamentally related by {membership}, and derivatively related by {application}}. Now, focus on the italicised sentence. To express what it needs to, our word ``membership'' must \emph{somehow} pick out the joint-carving, fundamental relation of membership, $\in$, rather than a derivative copy of that relation (such as $\in^{IJ}$). Indeed: \emph{somehow}. We must ask Seb: \emph{how}? How did it come to do this? Indeed, how could it \emph{so much as possibly} do this?

To be clear: I am not suggesting that we should start to worry that ``membership'' picks out some relation other than membership.\footnote{In particular, I want to cancel any suggestion that I am advocating for a model-theoretically inspired scepticism about meaning. I believe that such scepticism is deeply incoherent (see \cites[ch.7]{Button:LoR}{Button:BIVMT}[ch.9]{ButtonWalsh:PMT}).} My baffled and baffling questions for Seb arise out of his own commitments. Seb insists that membership \emph{is} joint-carving; and Seb must believe \emph{that} ``membership'' refers to membership; so Seb owes us and himself an explanation of \emph{why} ``membership'' refers to a joint-carving relation. 

Cutting to the chase: no good explanation is possible. Seb will end up saying that membership in particular, and joint-carving relations in general, are \emph{referentially magnetic}.\footnote{See e.g\ \textcite[23--33]{Sider:WBW}.} That last phrase has gained currency in contemporary philosophy, but we should be clear that Seb might as well have said that some relations \emph{cry out to be named};\footnote{Again, I am echoing Putnam's rhetorical criticism of magical theories of reference. I am also implicitly drawing upon arguments like  \textcites[ch.12]{Button:LoR}{Wrigley:SOII}.} that explanation stops here. Seb has been forced to posit a \emph{metasemantic surd}. Of course, there is no inconsistency in positing such a surd. But it is as absurd as anything I can imagine within mathematical metasemantics. 

In the end, Seb's world-view has so many degrees of freedom that it can only be held together by some (ab)surd(ity). And what goes for Seb goes for Fern. We must conclude that neither Seb nor Fern is right.

\subsection{Pluralistic joint-carving and its surds}\label{s:SameOntology}
This conclusion does not immediately entail that joint-carving realism  is wrong. Perhaps joint-carving realists can simply agree with us that Seb and Fern are both wrong. 

A terrible way to do this would be to claim that, in fact, $\Psi$s are fundamental, rather than sets or functions. Whatever exactly these $\Psi$s are claimed to be, this will amount to enriching our earlier dialogue with a third character, Sy, who offers $\Psi$-theoretic metaphysical foundations as a rival to set-theoretic and function-theoretic foundations.\footnote{Perhaps Sy endorses the Complemented Story, and \BLTzf, as presented in \textcite{Button:LT3}.} Barring some shocking logico-mathematical genius on Sy's part, this will be painfully epicyclic.

The only remaining option for the joint-carving realist is to insist that neither Seb nor Fern is right, because \emph{both} membership and application are joint-carving. On this view, set-theoretic and function-theoretic foundations present us with \emph{distinct} ontologies, but these ontologies are not \emph{competitors}.  Instead, we have two disjoint but perfectly fundamental hierarchies---a set-hierarchy and a function-hierarchy---sitting happily alongside each other in mathematical heaven. 

Making things personal again, let Pearl advocate this pluralistic version of joint-carving realism. In putting sets and functions on a par, Pearl has taken a step in the right direction. Alas, her continued adherence to joint-carving means that her position is still untenable, for both metaphysical and metasemantic reasons. 

On the \emph{metaphysical} front: Pearl's position is approximately as bad as Seb's and Fern's, but for slightly different reasons. Recall that Seb (and Fern) posited an inexplicable and for that reason unknowable fact: that membership rather than application (or vice versa) was fundamental. Pearl deliberately and explicitly avoids this, but in doing so she incurs the cost of profligacy at the fundamental level. After all, as we saw in \S\ref{s:judicial}, the synonymy result (my Main Theorem) shows that anything we can do with sets, we could instead do with functions (and vice versa), so there is no \emph{need} to treat both as fundamental; either alone would do.\footnote{Cf.\ the competition between \textsc{Austerity} and \textsc{Reasons} discussed in \textcite[\S8]{Button:SRSTP}.}

On the \emph{metasemantic} front: Pearl inherits all of Seb's and Fern's problems, unchanged. Exactly like Seb: Pearl will require that, \emph{somehow}, our word ``membership'' picks out the joint-carving membership relation, as found in the fundamental set-hierarchy, rather than (for example) the $J$-translation of that relation, as found in the fundamental function-hierarchy. And exactly like Fern: Pearl will require that, \emph{somehow}, our word ``application'' picks our the joint-carving application function , as found in the fundamental function-hierarchy, rather than (for example) the $I$-translation of that function, as found in the fundamental set-hierarchy. As before: my point is not that we should worry \emph{what} ``membership'' or ``application'' refer to, but that Pearl cannot possibly explain \emph{why} these words refer to (what she insists are) joint-carving stuff. She must posit a metasemantic surd; this is absurd.

Evidently, we must abandon Pearl's pluralism, just as we abandoned Seb's and Fern's monisms. And with this, joint-carving realism has run out of options. We must give up on the idea of joint-carving.

\section{Metaphysical foundations without joint-carving?}\label{s:metaphysicsnotdead}
Having given up on the idea of joint-carving, must we also abandon the very idea of providing metaphysical foundations for mathematics? Candidly: I don't know, but \emph{maybe} not. In this extremely speculative section, I will sketch an approach to metaphysical foundations which might just survive the demise of joint-carving realism. (For ease of readability, I will write as if I \emph{endorse} this approach; but I want to emphasize that I am just \emph{exploring} it.)

Suppose that, even after the discussion of \S\ref{s:metaphysics}, we still want set theory to supply us with metaphysical foundations. Given everything we have learned, we will have to accept that function theory also supplies us with metaphysical foundations. Moreover, we will eschew the idea that mathematical heaven has natural joints to carve. So we cannot say that some of the denizens of mathematical heaven are sets ``first'' (or ``fundamentally'') and functions ``second'' (or ``derivatively''). Instead, we will be forced to say that they are sets \emph{just as much} as they are functions. In slightly more detail, we will be forced to say something like this: 
\begin{quote}
	Mathematical heaven comprises a great many objects, which are richly structurable (see \eqref{n:things} of \S\ref{s:judicialfoundations}). Moreover, these objects are sets, arranged in an iterative hierarchy; and the fact that they are so arranged convinces us that our glimpses of mathematical are not mere mirages (see \eqref{n:consistent} of \S\ref{s:judicialfoundations}). Equally, though: these very same objects are functions, arranged in an iterative hierarchy. That is to say: the objects of mathematical heaven can be described and conceptualized equally well in at least two different (but synonymous) ways. They are sets, and they are functions, and they are not one ``before'' they are other.
\end{quote}
That is the approach. But it invites a puzzled question: \emph{How can anything be both a set and a function?} This question is important, and I will spend the rest of this section grappling with it.\footnote{Thanks to Chris Menzel and Chris Scambler for really pressing on this question.} 

We sometimes speak of sets and functions in cutesy terms: of sets as carrier bags full of objects (their members); of functions as things which, sea-cucumber-like, suck in arguments and spit out values. If we earnestly think about mathematical objects like this, then we will certainly be confused to be told that sets are also functions (carrier bags don't spit out objects). But, once we have explicitly spelled out this possible source of puzzlement, we can immediately set it aside. After all, such cutesy, carrier-bag-style imagery should never have been treated very seriously.\footnote{See \textcite[396--7]{OliverSmiley:CST}.}

We might instead think of sets and functions more \emph{structurally}. If we do that, it becomes totally unremarkable to be told that sets are also functions. Comparably: we do not even blink if we are told that someone is both a conductor and married to one of the orchestra's violinists. And note that this way of dissolving the puzzlement does not require that we think of sets and functions in \emph{wholly} structural terms: if we want, we can allow that sets and functions have some intrinsic features; we just need to accommodate their structural features when we characterize them \emph{as} sets or as functions.

We might also think of sets and functions \emph{presentationally}. Comparably: we would be perplexed to be told that something is both a duck and also a rabbit; but we are all familiar with \emph{duck--rabbits}, i.e.\ things which, viewed one way, are pictures of ducks, and viewed another, are pictures of rabbits. Inspired by duck--rabbits, we might say: mathematical heaven is a set--function-hierarchy: a bunch of entities which, viewed (or described) one way, are a hierarchy of sets, and viewed (or described) another, are a hierarchy of functions. 

The \emph{structural} and the \emph{presentational} approaches need not be competitors. Indeed, they are very naturally complementary. In drawing a duck--rabbit, you lay down a system of geometric shapes; the spatial relations between these shapes make the two viewpoints possible. Comparably: the structure of mathematical heaven makes possible its presentational multiplicities. 

Of course, the analogy between duck--rabbits and the (envisaged) set--function-hierarchy is limited in several ways, but one limitation merits particular comment. Presented with a duck--rabbit, we can describe it in \emph{duck}-terms (``its beak points left'') or in \emph{rabbit}-terms (``its ears point left''); but we can also describe it more \emph{neutrally}, exhaustively characterizing the underlying spatial structure by describing the various geometric shapes on the page (perhaps in terms of vector graphics). Now, we can certainly describe the set--function-hierarchy in \emph{set}-terms (via the Set Story), or in \emph{function}-terms (via the Function Story); but I very much doubt that the set--function-hierarchy can be adequately characterized in \emph{neutral} terms. 

This is not to say that there is anything \emph{inadequate} about the non-neutral descriptions that we might offer, i.e.\ our descriptions of the hierarchy in terms of sets or functions (or anything else). It is just to say that we cannot describe mathematical heaven without selecting some notational conventions. In turn: that is not (yet) to say that mathematical heaven is in any sense `constituted by' or `dependent upon' our conventional choices. It may just amount to the observation that our descriptions of mathematical heaven always bear our trace. 

I do not pretend that this is the last word on the question: \emph{How can anything be both a set and a function?} In the end, fully answering that question requires fully articulating this approach to metaphysical foundations. I have only offered a sketch of this approach here; I hope that this sketch is sufficient to show that the approach is worth pursuing.

\section{Conclusion, and Putnam revisited}\label{s:Putnam}
I began this paper with a quote from Putnam. That quote implicitly referred back to a comment he had offered nearly four decades earlier: 
\begin{quote}
	A mathematical example: you have a theory according to which mathematical heaven consists of objects called sets. The theory says that some of these sets are functions---in fact it says that those sets which are sets of ordered pairs satisfying a certain functionality condition are functions, and then it proceeds to translate your favorite textbook on calculus. The other theory says that mathematical heaven consists of things called functions,	even zero, $1$, $2$, $3$ turn out to be functions, and it says that some of these functions are sets---in fact those functions which take on only the values zero and $1$ are sets. Now, these are formally incompatible[. \ldots But they] are \emph{equivalent descriptions}[\ldots]. A sophisticated realist should not be bothered by the collapse of the `One True Theory' version of realism.\footnote{\textcite[287--8]{Putnam:CommentsOnKripke}.} 
\end{quote}
Putnam's remarks are suggestive, but somewhat opaque. What Putnam called the ``the `One True Theory' version of realism'' is just what I have called \emph{joint-carving realism}.\footnote{\textcites[xi, 49, 73, 143, 210]{Putnam:RTH}[352]{Putnam:MTFS} frequently described his opponents as believing that there is just `One True Theory', during his internal realist period (and cf.\ \cites[9, 201]{Button:LoR}). Note, though, the capital ``T'' on ``True'' (and cf.\ \cite[553--4]{Field:RR}). As flagged at the start of  \S\ref{s:metaphysics}, Sider's position is explicitly not that there is only one true theory, but that there is only one theory whose ``conceptual structure matches reality's structure''.} Moreover, Putnam held that joint-carving realism collapses, because set theory and function theory are ``\emph{equivalent descriptions}''. Unfortunately, he did not specify the relevant formal set theory, the relevant formal function theory, or the relevant notion of an \emph{equivalent description}.\footnote{Given his earlier \parencite*[7--9]{Putnam:MWF} use of the notion of ``equivalent description'', I conjecture that Putnam himself had in mind the (mere) \emph{mutual interpretability} of von Neumann's function theory and \Thnamed{NBG}. If so, then his line of thought was rather too quick, for reasons discussed in \S\ref{s:vonNeumann} (see also \cite[\S9.4, \emph{Second}]{Button:LT2}).} And whilst Putnam stated that a ``sophisticated realist'' could acknowledge that we are dealing with ``equivalent descriptions'', he did not tell us how, or what such sophisticated realism might amount to.\footnote{That quote, alongside other comments from the same conference \parencite*[esp.\ 228]{Putnam:RU}, were Putnam's first public avowal of (a version of) internal realism. But what \emph{that} amounted to was never straightforward; see \textcite[chs.8--11]{Button:LoR}.} 

This paper can be read as a reconstruction (and partial vindication) of Putnam's remarks. I introduced the bare idea of an iterative concept of \emph{function}, via the Function Story (see \S\ref{s:story}). By stages, I turned this Story into a quasi-categorical, categorial, formal theory, \FLT. And an extension of this theory has all the power of \ZF; indeed, \ZF and \FLTzf are synonymous (see \S\ref{s:flt}).

This perfectly precise equivalence shows that set theory and function theory provide exactly the same judicial foundations for mathematics, subject only to a change in notation (see \S\ref{s:judicial}). Moreover, when this point is fully appreciated, it leads to the collapse of joint-carving metaphysical foundations (see \S\ref{s:metaphysics}). But we may yet be left with the possibility of a ``sophisticated realism'', according to which set theory and function theory present us with the very same metaphysical foundation (see \S\ref{s:metaphysicsnotdead}).

\section*{Acknowledgements}
In late 2019, I developed both \FST and \FLT, and sketched the proof of my Main Theorem (see \S\ref{s:de:sketch}). I sent this to Thomas Forster and Albert Visser at the time, and I would like to thank them both for early discussion and inspiration: the paper was originally motivated by their suggestion that I should try to offer a level-theoretic axiomatization of multisets. (This is indeed possible; see \cite{Button:WST}.) 

In June 2021, I presented this paper at the seminar series which gave rise to this volume, covering the formal theories, the proof-strategy for the Main Theorem, and the significance of all this for joint-carving realism. I wish to thank that audience for helpful comments, but most especially Neil Barton and Giorgio Venturi. Huge thanks also to Chris Menzel, and to an anonymous referee for this volume.

\startappendix

\section{Background on \LT}\label{s:app:lt}
The remainder of this paper consists of appendices concerning \FLT. Recall that \FLT axiomatizes the Function Story of \S\ref{s:story}. Technically, it is a straightforward adaptation of the theory \LT, which axiomatizes the Set Story. I draw on \LT in my discussion of \FLT and in several footnotes in the main text. So, in this appendix, I will quickly introduce \LT.\footnote{For full discussion of \LT, see \textcite{Button:LT1}.} I begin with the key definition (the paranthetical ``\ref{ext}'' again indicates that this definition relies upon \ref{ext}.) 
\begin{define}[\ref{ext}]\label{def:lt:level}
	For any $a$, let $\pot{a} = \Setabs{x}{(\exists c \in a)x  \subseteq c}$, if it exists.
	
	Say that $h$ is a \emph{history}, written $\histpred(h)$, iff $(\forall x \in h)x = \pot{(x \cap h)}$.
	
	Say that $s$ is a \emph{level}, written $\levpred(s)$, iff $\exists h(\histpred(h) \land s = \pot{h})$. \defendhere
\end{define}\noindent
Essentially, the levels serve as surrogates for the stages in the Set Story. We can now set up some axioms: 
\begin{listaxiom}
	\labitem{Extensionality}{ext} $\forall a \forall b \big(\forall x(x \in a \liff x \in b) \lonlyif a =b\big)$
	\labitem{Separation}{sep} $\forall F \forall a(\Setabs{x \in a}{F(x)}\text{ exists})$ \hfill \emph{with ``$F$'' a predicate-variable}
	\labitem{Stratification}{lt:strat} $\forall a\exists s\big(\levpred(s) \land a \subseteq s\big)$
	\labitem{Endless}{lt:cre} $\forall s \exists t \ s \in t$
	\labitem{Infinity}{lt:inf} $(\exists s \neq \emptyset)(\forall q\in s)\exists r(q \in r \in s)$
	\labitem{Unbounded}{lt:rep} $\forall P\forall a \exists s(\forall x \in a) P(x) \in s$ \hfill \emph{with ``$P$'' a function-variable}
\end{listaxiom}
Now \LT itself is just \ref{ext} + \ref{sep} + \ref{lt:strat}. Then \LTplus is \LT + \ref{lt:cre}, and \ZF is equivalent to \LT + \ref{lt:inf} + \ref{lt:rep}. 

Using just \ref{ext} and \ref{sep}, we can prove that the levels are well-ordered by membership. This licenses the introduction of a powerful tool: 
\begin{define}[\LT]\label{def:lt:levof}
	Let $\levof{a}$ be the $\in$-least level including $a$. So, $a \subseteq \levof{a}$ and $\lnot (\exists s \in \levof{a})(\levpred(s) \land a \subseteq s)$. \defendhere
\end{define}\noindent
Note that $\levof{a}$ obeys various intuitive principles, e.g.\ that if $a \in b$ then $\levof{a} \in \levof{b}$.\footnote{\label{fn:levwellbehaved}See \textcite[Lemma 3.12]{Button:LT1}.} 

\section{Elementary considerations regarding \FLT}\label{app:FLT}
I now turn to considering \FLT. I stated its axioms in \S\ref{s:flt:formal}, but there are a few things I need to clear up. 

First: I need to be a little more precise about the background logic for handling partial functions (see \S\ref{s:pedantry:partial}). Second: in \S\ref{s:flt:formal}, I invoked the predicate, $\fevpred$, i.e.\ the formal notion of a fevel (i.e.\ a functional-level); but I still need to define it (see \S\ref{s:elementary:definitions}). 

\subsection{Formalizing partial functions}\label{s:pedantry:partial}
There is a ``pseudo-partial'' approach to partial functions, according to which, to say that $f$ is undefined on input $x$ is to say that $f(x)=\bot$, where $\bot$ is some ``default object''. I say that this is \emph{pseudo}-partial, since when $\bot$ is counted as a value, this makes all functions total. I have no interest in this approach. I want to present a theory of genuinely partial functions. I really want to allow that $f$ may be undefined on input $x$; that $\forall y\ f(x) \neq y$.

An obvious way to implement this is via a negative free logic (with no non-logical primitives). In brief, the syntax for this approach will be as follows:
\begin{listbullet}
	\item Variables (and nothing else) are primitive terms; $\tau(\sigma)$ is a term whenever both $\tau$ and $\sigma$ are; nothing else is a term. 
	\item For any term $\tau$, we define $\isdefined{\tau}$ via $\tau=\tau$. 
	\item  $\isdefined{\text{x}}$ for any variable $\text{x}$; if $\alpha$ is an atomic formula containing a term $\tau$, then $\alpha \lonlyif \isdefined{\tau}$; and $\isdefined{\tau(\sigma)} \rightarrow (\isdefined{\tau} \land \isdefined{\sigma})$ for any terms $\tau$ and $\sigma$. 
\end{listbullet}
Implicitly, this is how I formulated \FLT (and its extensions) in the main paper. Once you get used to it, it is quite neat.

There is a downside to this approach. Whilst ``$\tau(\sigma)$'' intuitively stands for the application of $\tau$ to $\sigma$, no object-language symbol expresses \emph{application}. This makes it hard to consider interpretations of theories offered with this formalism. And, of course, the main focus of this paper is less on \FLT itself, but more on the synonymy between \FLTzf and \ZF, which requires considering interpretations. To remedy this, we can simply add a primitive two-place function symbol, ``$\appfunction$'', to express application. We would now regard ``$\tau(\sigma)$'' as an abbreviation of ``$\appof{\tau}{\sigma}$'', which is a term whenever both $\tau$ and $\sigma$ are. The result is a negative free logic with a primitive function-symbol. 

Still, when our aim is to prove synonymy, it is easiest to work with a theory in classical logic (rather than free logic) which uses only primitive \emph{relation} symbols (instead of \emph{function} symbols). With that in mind, instead of expressing application via a primitive function symbol, we can express it with a primitive three-place predicate, ``\valfunction''. We then need an axiom stating that \valfunction is functional:
\begin{listaxiom}
	\labitem{FunVal}{flt:valfun} $\forall f \forall x\forall y \forall z\big((\valfunction(f,x,y) \land \valfunction(f,x,z))\lonlyif y = z\big)$
\end{listaxiom}
This allows us to abandon free logic.\footnote{There is a subtlety: free logics usually allow for empty domains. In abandoning free logic, we rule out the empty domain. This wrinkle is, however, insignificant. Those who find it causes issues can remedy it in either of two ways: (1) add an axiom $\exists x\ x = x$ to \FLT (when formulated functionally); (2) when considering the interpretation $J : \LTplus \functionto\FLTplus$, (re)formulate \LT in a free logic.} We can easily map from the functional approach to the relational approach, by (for example) treating ``$f(x) = y$'' as ``$\valfunction(f, x, y)$'', and ``$f(x) \simeq g(x)$'' as ``$\forall y(\valfunction(f,x,y) \liff \valfunction(g,x,y))$''. 

My attitude is: all these formalisms are equivalent; use whichever formalism is easiest for the purposes at hand. 

\subsection{Key definitions and axioms}\label{s:elementary:definitions}
My next task is to define the key predicate, $\fevpred$, i.e.\ the formal notion of a fevel, which we use in \FLT's axiom \ref{flt:strat}. 

Recall from Definition \ref{def:funset} that \funsets are partial identity functions. It follows from \ref{flt:ext} that \funsets obey extensionality$_{\funin}$.\footnote{i.e.\ $\forall x(x \funin f \liff x \funin g) \lonlyif f = g$, for any \funsets $f$ and $g$.} In turn, this licenses us in introducing a function-theoretic analogue of ordinary set-builder notation:
\begin{define}[\ref{flt:ext}]
	If it exists, then $\Unitabs{x}{\phi}$ is the \funset given as follows:
	\begin{align*}
		\text{if }\suballforin{\phi}{v}{x}\text{, then }&\Unitabs{x}{\phi}(v)=v\\
		\text{otherwise, }&\Unitabs{x}{\phi}(v)\text{ is undefined}
	\end{align*}
	I tweak this notation in obvious ways, e.g.\ allowing $\Unitabs{x \funin a}{\phi} = \Unitabs{x}{x \funin a \land \phi}$. \defendhere
\end{define}\noindent
Using this \funset-builder notation, here is the define of a \emph{fevel}: 
\begin{define}[\ref{flt:ext}]\label{def:flt:fevel}
	For any $f$, let $\hfpot{f} = \Unitabs{x}{(\exists g \funin f)x \funsubeq g}$, if it exists.\footnote{Note that we do not initially assume that $\hfpot{f}$ exists, for all $f$. }
	
	Say that $h$ is a \emph{functional-history}, written $\fistpred(h)$, iff $(\forall x \funin h)x = \hfpot\Unitabs{z \funin h}{z \funin x}$.
	
	Say that $s$ is a \emph{fevel}, written $\fevpred(s)$, iff $\exists h(\fistpred(h) \land s = \hfpot{h})$. \defendhere
\end{define}\noindent
The utility of this definition is far from immediately obvious! But a quick glance will confirm that it arises from Definition \ref{def:lt:level}, by translating set-theoretic into \funset-theoretic vocabulary. 

\subsection{Key results}\label{s:elementary:results}
Of course, we need to check that Definition \ref{def:flt:fevel} does what we want. Ultimately, it is vindicated by Theorem \ref{thm:fst-flt}, which tells us that \FLT and \FST make exactly the same claims about functions. I will not prove Theorem \ref{thm:fst-flt} here;\footnote{This is proved as for \LT; see \textcite[\S4]{Button:LT1}.} but I will note that, on the way to that result, we will want to prove the fundamental theorem of fevel theory:
\begin{thm*}[\textbf{\ref{thm:flt:wo}}; \ref{flt:ext}, \ref{flt:comp}] The fevels are well-ordered by $\funin$.
\end{thm*}
\begin{proof}[Proof sketch]
	The first step is to show that \funsets obey Separation$_{\funin}$. So we must show that $\Unitabs{x \funin f}{F(x)}$ exists for any property $F$ and any function $f$. To do this, fix a second-order map $P$ such that $P(x) = x$ if $x \funin f \land F(x)$, and $P(x)$ is undefined otherwise; now use \ref{flt:comp} and \ref{flt:ext}.
	
	From here, the remaining steps are as for \LT (i.e.\ see \cite[\S3]{Button:LT1}). We first show that the fevels are well-founded. This allows us to prove that any member$_{\funin}$ of a fistory is a fevel. We can then show that a fevel is precisely the result of hitting the \funset of all earlier fevels with $\hfpot$, i.e.\ that $\fevpred(f)$ iff $f = \hfpot{\Unitabs{s \funin f}{\fevpred(s)}}$. Using this and induction on $\funin$ we can show that fevels are linearly ordered by $\funin$.
\end{proof}\noindent
The well-ordering of the fevels also allows us to introduce a powerful notion (compare Definition \ref{def:lt:levof}): 
\begin{define}[\FLT]\label{def:fevof}
	Let $\fevof{f}$ be the $\funin$-least fevel including $f$. So $f \funsubeq \fevof{f}$ and $\lnot (\exists s \funin \fevof{f})(\fevpred(f) \land f \funsubeq s)$. \defendhere
\end{define}\noindent
Together with Theorem \ref{thm:flt:wo}, this licenses $\funin$-induction. Moreover, the fevels of \FLT, like the levels of \LT, obey some very intuitive principles, e.g.\ that if $g \funin f$, then $\fevof{g} \funin \fevof{f}$.\footnote{It is easy to obtain a result like \textcite[Lemma 3.12]{Button:LT1}.} 

\subsection{\FLT and categories}\label{app:FLT:cat}
The final key result I mentioned, concerning \FLT, is that any cumulative hierarchy of functions is a \emph{category}. To show this, it will help to introduce some abbreviations.
\begin{define}[\ref{flt:ext}]\label{def:composition}
	Let $\dom(f) = \Unitabs{x}{\exists v\ f(x) = v}$, and $\cod(f) = \Unitabs{x}{\exists v\ f(v) = x}$. Say that $f$ is an \emph{arrow} from $\dom(f)$ to $\cod(f)$; we may write $f : a \functionto b$ to indicate that $\dom(f) = a$ and $\cod(f) = b$. When $f : a \functionto b$ and $g : b \functionto c$, the \emph{composite} of $f$ and $g$ is the function $g \circ f = \Lamabs{x}{g(f(x))}: a \functionto c$. \defendhere
\end{define}
\begin{lem}[\FLT]\label{lem:flt:category}
	For any $f, g, h$:
	\begin{listn-0}
		\item\label{cat:domcod} $\dom(f)$ and $\cod(f)$ exist
		\item\label{cat:compexists} If $\cod(f) = \dom(g)$, then $g \circ f$ exists
		\item\label{cat:funset} $f \circ a = f = b \circ f$, when $f : a \functionto b$ 
		\item\label{cat:ass} $h \circ (g \circ f) = (h \circ g) \circ f$, when $f : a \functionto b, g : b \functionto c, h : c \functionto d$.  
	\end{listn-0}
\end{lem}
\begin{proof}
	\emphref{cat:domcod} By Separation$_{\funin}$ on $\fevof{f}$. 
	
	\emphref{cat:compexists} Let $P(x)$ be given by $g(f(x))$; use \ref{flt:comp} on the greater of $\fevof{f}$ and $\fevof{g}$. 
	
	\emphref{cat:funset} Using \eqref{cat:compexists}, note that $f \circ a$ exists, since $\cod(a) = a = \dom(f)$ by \ref{flt:ext}. Now $(f \circ a)(x) \simeq f(a(x)) \simeq f((\dom(f))(x)) \simeq f(x)$ for all $x$; so $f \circ a = f$ by \ref{flt:ext}. The case for $b$ is similar.  
	
	\emphref{cat:ass} 
	By \ref{flt:ext}.
\end{proof}\noindent
This lemma immediately yields Theorem \ref{thm:flt:category}. However, note that Definition \ref{def:composition} in effect identifies a function's range with its codomain. And this has certain immediate limiting effects. For example: the category we obtain from \FLTplus is \emph{not} equipped with any of: an initial object; a terminal object; products; equalizers; or pullbacks.

The case of products is worth commenting on in a little detail. Having decided how to implement ordered pairs (see \S\ref{app:Judicious}), we can take a cue from set theory, and define $f \times g \coloneq \Unitabs{\tuple{x,y}}{x \funin f \land y \funin g}$. So defined, \FLTplus proves that this exists for any $f$ and $g$. 
	%
	%
But this is almost never a product in the \emph{category theorist's} sense.\footnote{See e.g.\ \textcite[38--41]{Awodey:CT}.} To see why, note the following.
%
\begin{lem}[\FLTplus]
	Let $A$ and $B$ be \funsets with a product diagram $A \xleftarrow{\pi_1} P \xrightarrow{\pi_2} B$. Then $\max(|A|, |B|) = |A \times B|$,\footnote{Here I write $|Z|$ for $Z$'s cardinality, and ``$\times$'' is as defined a few lines ago.} so that either $|A| < 2$ or $|B| < 2$.
\end{lem}
\begin{proof}
	By the UMP for products, for any $A \xleftarrow{q_1} Q \xrightarrow{q_2} B$, there is a unique $u : Q \functionto P$ such that $q_1 = \pi_1 \circ u $ and $q_2 = \pi_2 \circ u$. And since ranges are codomains, $|P| \leq |Q|$. 
		
	I will show that $|P| = \max(|A|, |B|)$. Since ranges are codomains, for any \funset $Q$ there are arrows $q_1 : Q \functionto A$ and $q_2 : Q \functionto B$ iff $\max(|A|, |B|) \leq |Q|$. Since $\pi_1$ and $\pi_2$ exist, $\max(|A|, |B|) \leq |P|$. Now fix some diagram $A \xleftarrow{q_1} Q \xrightarrow{q_2} B$ with $|Q| = \max(|A|, |B|)$, and observe that  $\max(|A|, |B|) \leq |P| \leq |Q| = \max(|A|,|B|)$.
	
	I will now show that $|P| = |A \times B|$. Fix the diagram $A \xleftarrow{p_1} A \times B \xrightarrow{p_2} B$ via $p_1(\tuple{x,y}) = x$ and $p_2(\tuple{x,y}) = y$. Obtain $u : A \times B \functionto P$ by the UMP. Suppose that $u(\tuple{a,b}) = u(\tuple{a',b'})$; then $a = p_1(\tuple{a,b}) = \pi_1(u(\tuple{a,b})) = \pi_1(u(\tuple{a',b'}))= p_1(\tuple{a',b'}) = a'$; and similarly $b = b'$; so $\tuple{a,b} = \tuple{a',b'}$. Generalizing, $u$ is an injection, so $|A \times B| \leq |P| \leq |A \times B|$.
\end{proof}\noindent
So, although \FLTplus merrily gives us \funsets of ordered pairs, it gives us almost no categorial products. 

If we want to have nice categorial things, like products and equalizers, we may want to move away from \FLT and its extensions. In turn, this requires that we modify the Function Story (from \S\ref{s:story}). Specifically, we could tweak the Story to say that, at any stage, we find all functions whose domains and \emph{codomains} are exhausted by the functions we found earlier. This will allow that codomains can outstrip ranges, and thereby (ultimately) deliver a topos. 

\section{Autonomous judicial foundations in \FLT}\label{app:Judicious}
In \S\ref{s:foundationsketch}, I noted that we can provide ``autonomous'' judicial foundations within \FLT. Essentially, we can take inspiration from \FLT's category-theoretic nature, to develop function-theoretic implementations of various mathematical notions (whose set-theoretic implementations are perfectly familiar). Here are some ways one might proceed.

\emph{Arbitrary functions.} In \FLT, we take $1$-place functions as our basic currency. But it is worth considering $n$-place functions. Set theorists usually regard $n$-place functions as sets of $n\mathord{+}1$-tuples satisfying a functionality constraint. Function-theorists will probably prefer to curry. For example, to handle 3-place functions, we might say: 
\begin{align*}
	f(x, y, z)\text{ abbreviates }((f(x))(y))(z)
\end{align*}
Generalizing, function-theoretic foundations allow us to regard \emph{any} object (i.e.\ any one-place function) immediately as a (partial) $n$-place function, for any $n$.

\emph{Arbitrary relations.} Set theorists usually regard $n$-place relations as sets of $n$-tuples. Taking the above hint, function-theorists can immediately regard any object as an $n$-place relation, for any $n$. For example, to handle 3-place relations, we might say: 
\begin{align*}
	\emph{Rel}_f[x,y,z]\text{ abbreviates }(\exists g \funin f)g(x, y)=z
\end{align*}
where the right-hand-side uses currying. But we can now read ``$\emph{Rel}_f[x,y,z]$'' as ``the relation (associated with) $f$ holds of $x$, $y$ and $z$ (in that order)''. 

\emph{Ordered pairs.} Set theorists usually render ordered pairs via Kurotawksi's fairly arbitrary definition that $\tuple{a, b} \coloneq \{\{a\}, \{a,b\}\}$. Function theorists might instead (fairly arbitrarily) say that $\tuple{a, b}$ is the function which maps $a \mapsto b \mapsto b$, and is undefined otherwise. (However, ordered pairs may not play such a dominant role in function-theoretic foundations, given the earlier remarks about arbitrary functions/relations.)

\emph{Ordinals.} Set theorists usually adopt von Neumann's implementation of the ordinals. Function-theorists might prefer to implement ordinals via the following  intuitive thought: each ordinal $\beta$ is the function which maps $\alpha \mapsto \alpha \dotminus 1$ for all $\alpha < \beta$, and which is undefined otherwise. (Here I am using $\dotminus$ for truncated subtraction; this can of course be given a more rigorous implementation.) 

These are just a few examples of mathematical concepts which can be given autonomous, function-theoretic implementations. We could continue, but I leave this for anyone who really wants to put \FLT to work.

\stopappendix

\printbibliography

@Article{Button:WST,
	Title = {Wand/Set Theory},
	Author = {Tim Button},
	Journal = {arXiv:2308.06789},
	Year = {2023}}

@InCollection{Button:BIVMT,
	Title = {Brains in vats and model theory},
	Author = {Tim Button},
	Booktitle = {The brain in a vat},
	Editor = {Sanford Goldberg},
	Publisher = {Cambridge University Press},
	Address = {Cambridge},
	Pages = {131--54},
	Year = {2016}}

@InCollection{Bouvere:ST,
	Title = {Synonymous Theories},
	Year = {1963},
	Booktitle = {The Theory of Models: Proceedings of the 1963 International Symposium at {Berkeley}},
	Author = {Karel de Bouv\`ere},
	Publisher = {North-Holland},
	Address = {Amsterdam},
	Editor = {J W {Addison}, Leon {Henkin} and Alfred {Tarski}}}

@PhDThesis{Montague:PhD,
	Title = {Contributions to the axiomatic foundations of set theory},
	Year = {1957},
	School = {Berkeley},
	AUthor = {Montague, Richard M}}

@Article{Santanna:Flow,
	Author = {Adonai S {Sant'Anna} and Ot\'avio {Bueno} and Marcio P P {de França} and Renato {Brodzinski}},
	Title = {Flow: The Axiom of Choice is independent from the Partition Principle},
	Journal = {arXiv:2010.03664},
	Year = {2020}}

@Article{OliverSmiley:CST,
	Author = {Alex {Oliver} and Timothy {Smiley}},
	Title = {Cantorian Set Theory},
	Journal = {Bulletin of Symbolic Logic},
	Year = {2018},
	Volume = {24},
	Number = {4},
	Pages = {393--451}}

@InCollection{Amijee:PSR,
	Author = {Amijee, Fatema},
	Title = {Principle of Sufficient Reason},
	Booktitle = {The Routledge Handbook of Metaphysical Grounding},
	Address = {New York},
	Publisher = {Routledge},
	Pages = {63--75},
	Year = {2020}}

@Article{Button:SRSTP,
	Author = {Button, Tim},
	Title = {Symmetric relations, symmetric theories, and {Pythagrapheanism}},
	Journal = {Philosophy and Phenomenological Research},
	Year = {2023},
	Pages = {583--612},
	Volume = {107},
	Number = {3}}

@Book{Awodey:CT,
	Author = {Awodey, Steve},
	Title = {Category Theory},
	Year = {2010},
	Edition = {2},
	Publisher = {Oxford University Press},
	Address = {Oxford}}

@Article{Robinson:TCMVNS,
	Author = {Robinson, Raphael M},
	Title = {The Theory of Classes: A Modification of von {Neumann}'s System},
	Journal = {Journal of Symbolic Logic},
	Year = {1937},
	Volume = {2},
	Number = {1},
	Pages = {29--36}}

@Article{Maddy:STF,
	Author = {Maddy, Penelope},
	Title = {Set-Theoretic Foundations},
	Journal = {Contemporary Mathematics},
	Year = {2017},
	Volume = {690},
	Pages = {289--322}}

@Book{Lewis:OPW,
	title         = {On the Plurality of Worlds},
	publisher     = {Basil Blackwell},
	year          = {1986},
	author        = {Lewis, David},
	address       = {Oxford},
}

@Article{Quine:OI,
	Author = {Quine, Willard van Orman},
	Title = {Ontology and Ideology},
	Journal = {Philosophical Studies},
	Year = {1951},
	Volume = {2},
	Number = {1},
	Pages = {11--15}}

@InCollection{Button:MIR,
	Author = {Tim Button},
	Title = {Mathematical internal realism},
	Booktitle = {Engaging Putnam},
	Editor = {James {Conant} and Sanjit {Chakraborty}},
	Publisher = {De Gruyter},
	Year = {2022},
	Pages = {157--82},
	Address = {Berlin}}

@InCollection{Isaacson:RMCST,
	Author = {Daniel {Isaacson}},
	Title = {The reality of mathematics and the case of set theory},
	Year = {2011},
	Editor = {Zsolt {Nov\`ak} and Andr\`as {Simonyi}},
	Booktitle = {Truth, Reference and Realism},
	Publisher = {CEU Press},
	Address = {New York}}

@Article{Field:RR,
	Author = {Hartry {Field}},
	Title = {Realism and Relativism},
	Journal = {The Journal of Philosophy},
	Year = {1982},
	Volume = {79},
	Number = {10},
	Pages = {553--67}}

@Book{Putnam:RTH,
	Author = {Hilary Putnam},
	Title = {Reason, Truth, and History},
	Year = {1981},
	Publisher = {Cambridge University Press},
	Address = {Cambridge}}

@InCollection{Putnam:MTFS,
	Author = {Hilary Putnam},
	Title = {Model Theory and the ``Factuality'' of Semantics},
	Booktitle = {Reflections on {Chomsky}},
	Editor = {Alexander George},
	Publisher = {Blackwell},
	Year = {1989},
	Address = {Oxford},
	Pages = {213--32}
}

@Book{Shapiro:PM,
	Author = {Stewart Shapiro},
	Title = {Philosophy of Mathematics: Structure and Ontology},
	Year = {1997},
	Publisher = {Oxford University Press},
	Address = {Oxford}
}

@Book{Button:LoR,
	Author = {Tim Button},
	Title = {The Limits of Realism},
	Year = {2013},
	Publisher = {Oxford University Press},
	Address = {Oxford}}

@Book{Wrigley:SOII,
	Author = {Wesley Wrigley},
	Title = {{Sider}'s {Ontologese} introduction instructions},
	Year = {2018},
	Journal = {Theoria},
	Pages = {295--308},
	Volume = {84}, 
	Issue = {4}}

@Article{Maddy:BA1,
	Author = {Penelope Maddy},
	Title = {Believing the Axioms {I}},
	Journal = {Journal of Symbolic Logic},
	Year = {1988},
	Volume = {53},
	Number = {2},
	Pages = {481--511}}

@Article{Maddy:BA2,
	Author = {Penelope Maddy},
	Title = {Believing the Axioms {II}},
	Journal = {Journal of Symbolic Logic},
	Year = {1988},
	Volume = {53},
	Number = {3},
	Pages = {736--64}}

@Book{Maddy:NM,
	Author = {Penelope Maddy},
	Title = {Naturalism in Mathematics},
	Year = {1997},
	Publisher = {Oxford University Press},
	Address = {Oxford}}

@Article{Sider:WBW,
	Author = {Sider, Theodore},
	Title = {Writing the Book of the World},
	Year = {2011},
	Publisher = {Oxford University Press},
	Address = {Oxford}}

@InCollection{Putnam:RU,
	Author = {Putnam, Hilary},
	Title = {Reference and Understanding, and Reply to {Dummett}'s Comments},
	Pages = {199--217, 226--8},
	Booktitle = {\emph{\textcite{Margalit:MU}}},
	Year = {1979}}

@InCollection{Putnam:CommentsOnKripke,
	Author = {Putnam, Hilary},
	Title = {Comments (on {Kripke}'s ``{A} puzzle about belief'')},
	Pages = {284--8},
	Booktitle = {\emph{\textcite{Margalit:MU}}},
	Year = {1979}}

@Article{vonNeumann:EAM,
	Author = {John {von Neumann}},
	Title = {Eine {Axiomatisierung} der {Mengenlehre}},
	Journal = {Journal f\"ur reine und angewandte Mathematik},
	Volume = {154},
	Year = {1925},
	Pages = {219--240},
	Notes = {Translated in Jean van Heijenoort, \emph{From Frege to G\"odel}, 1967, Harvard University Press (secondary page references to the English translation, which I depart from occasionally)}}

@Article{vonNeumann:DAM,
	Author = {John {von Neumann}},
	Title = {Die {Axiomatisierung} der {Mengenlehre}},
	Journal = {Mathematische Zeitschrift},
	Volume = {27},
	Year = {1928},
	Pages = {669--752},
}

@Article{vonNeumann:WAM,
	Author = {John {von Neumann}},
	Title = {\"Uber eine {Widerspruchsfreiheitsfrage} in der axiomatischen {Mengenlehre}},
	Journal = {Journal f\"ur reine und angewandte Mathematik},
	Year = {1929},
	Volume = {160},
	Pages = {227--41}
}

@Article{FriedmanVisser:WBIS,
	Author = {Harvey M {Friedman} and Albert {Visser}},
	Title = {When bi-interpretability implies synonymy},
	Journal = {Logic Group Preprint Series}, 
	Volume = {320},
	Pages = {1--19}, 
	Year = {2014}}

@online{Putnam:blog,
	Author = {Putnam, Hilary},
	url = {http://putnamphil.blogspot.com/2014/12/},
	title = {\emph{Three blog posts: `The modal logical interpretation and ``equivalent descriptions''{}' (11.Dec.2014); `Continuing' (12.Dec.2014); and `Mathematical ``existence''{}' (13.Dec.2014)}},
	year = {2014}
}

@Article{Putnam:MWF,
	Author = {Putnam, Hilary},
	Title = {Mathematics without Foundations},
	Journal = {Journal of Philosophy},
	Year = {1967},
	Volume = {64},
	Number = {1},
	Pages = {5--22}}

@Book{ButtonWalsh:PMT,
  title     = {Philosophy and Model Theory},
  publisher = {Oxford University Press},
  year      = {2018},
  author    = {Tim {Button} and Sean {Walsh}},
  address   = {Oxford},
}

@Article{Boolos:ICS,
  author  = {Boolos, George},
  title   = {The Iterative Conception of Set},
  journal = {The Journal of Philosophy},
  year    = {1971},
  volume  = {68},
  number  = {8},
  pages   = {215--31},
}

@InCollection{Scott:AST,
  author    = {Scott, Dana},
  title     = {Axiomatizing Set Theory},
  booktitle = {Axiomatic Set Theory II},
  publisher = {American Mathematical Society},
  year      = {1974},
  editor    = {Thomas Jech},
  pages     = {207--14},
  note      = {Proceedings of the Symposium in Pure Mathematics of the American Mathematical Society, July--August 1967},
}

@InCollection{Shoenfield:AST,
  author    = {Shoenfield, Joseph R},
  title     = {Axioms of set theory},
  booktitle = {Handbook of Mathematical Logic},
  publisher = {North-Holland},
  year      = {1977},
  editor    = {Barwise, Jon},
  pages     = {321--44},
  address   = {London},
}

@Article{Button:LT1,
  author = {Button, Tim},
  title  = {Level Theory, Part 1: Axiomatizing the bare idea of a cumulative hierarchy of sets},
  year   = {2021},
  Volume = {27},
  Number = {4}, 
  Pages = {436--60},
  Journal = {Bulletin of Symbolic Logic}
}

@Article{Button:LT2,
	author = {Button, Tim},
	title  = {Level Theory, Part 2: Axiomatizing the bare idea of a potential hierarchy},
	year   = {2021},
	Journal = {Bulletin of Symbolic Logic},
	Volume = {27},
	Number = {4},
	Pages = {461--84}
}

@Article{Button:LT3,
	author = {Button, Tim},
	title  = {Level Theory, Part 3: A boolean algebra of sets arranged in well-ordered levels},
	year   = {2022},
	Journal = {Bulletin of Symbolic Logic},
	Volume = {28},
	Number = {1},
	Pages = {1--26}}
\end{document}